\documentclass{amsart}
\usepackage{amssymb,amsmath}

\newtheorem{theorem}{Theorem}[section]
\newtheorem{proposition}[theorem]{Proposition}
\newtheorem{lemma}[theorem]{Lemma}
\newtheorem{corollary}[theorem]{Corollary}

\theoremstyle{definition}
\newtheorem{definition}[theorem]{Definition}
\newtheorem{example}[theorem]{Example}

\theoremstyle{remark}
\newtheorem{remark}[theorem]{Remark}
\numberwithin{equation}{section}
\def\abs#1{\left| #1 \right|}
\def\cl{c\ell}
\def\coneg#1{\mathrm{C_G}(\mathbb{B},(#1))}
\def\C{\mathbb{C}}
\def\dbar{d\hspace*{-0.08em}\bar{}\hspace*{0.1em}}
\def\diff{\mathrm{Diff}}
\def\diffC#1{\diff^{#1}(\mathbb{B})_{\mathrm{cone}}}
\def\edgeg#1#2{\mathrm{R}^{#1}_{\mathrm{G}}(\gL,(#2))}
\def\emb{\hookrightarrow}
\def\eps{\varepsilon}
\def\etA{e^{-tA}}
\def\etgl{e^{-t\gl}}
\def\g{\mathfrak{r}}
\def\gG{\Gamma}
\def\gL{\Lambda}
\def\gLint{\overset{{}_\circ}{\gL}}
\def\gO{\Omega}
\def\gU{\Upsilon}
\def\ga{\alpha}
\def\gb{\beta}
\def\gd{\delta}
\def\gg{\gamma}
\def\gk{\kappa}
\def\gl{\lambda}
\def\go{\omega}
\def\gp{\varphi}
\def\gr{\varrho}
\def\gs{\sigma}
\def\gt{\theta}
\def\hsg{\mathcal{H}^{s,\gamma}}
\def\htensor{\hat{\otimes}_{\scriptscriptstyle{H}}}
\def\iBB{\mathrm{int}\,\mathbb{B}}
\def\intU{\int_{\gU}\!}
\def\inttU{\int_{t^{-1}\gU}\!}
\def\ladm{[\gl]_m}
\def\norm#1{\left\| #1 \right\|}
\def\N{\mathbb{N}}
\def\onull{\setminus 0}
\def\op{\mathrm{op}}
\def\pitensor{{\hat{\otimes}_\pi}}
\def\rpbar{\overline{\R}_+}
\def\R{\mathbb{R}}
\def\sK{\gs_{\wedge}}
\def\sM{\gs_{M}}
\def\sod{\,\boldsymbol{|}\ }
\def\spec{\operatorname{spec}}
\def\spk#1{\left<#1\right>}
\def\sumj#1{\sum_{j=0}^{#1}}
\def\sumk#1{\sum_{k=0}^{#1}}
\def\supp{\operatorname{supp}}
\def\tr{\operatorname{Tr}}

\begin{document}
\title[Full expansion of the heat trace for cone operators]{Full asymptotic
expansion of the heat trace for\\ non-self-adjoint elliptic cone operators}
\author{Juan B. Gil}
\address{Department of Mathematics\\ Temple University\\ 
Philadelphia, PA 19122, USA}

\begin{abstract}
The operator $\etA$ and its trace $\tr\etA$, for $t>0$, are
investigated in the case when $A$ is an elliptic differential operator on a
manifold with conical singularities. Under a certain spectral condition
(parameter-ellipticity) we obtain a full asymptotic expansion in $t$ of
the heat trace as $t\to 0^+$. As in the smooth compact case, the problem is
reduced to the investigation of the resolvent $(A-\lambda)^{-1}$. The main
step consists in approximating this family by a parametrix of
$A-\lambda$ constructed within a suitable parameter-dependent calculus. 
\end{abstract}

\maketitle

\section{Introduction}
In this paper the operator $\etA$, $t>0$, is investigated on manifolds with
conical singularities. The operator $A$ is assumed to be an elliptic 
differential operator of arbitrary positive order, not necessarily
self-adjoint, but satisfying an analog of Agmon's condition
(parameter-ellipticity) on a sector 
$\{\gl\in\C\sod \gp\le \arg\lambda\le 2\pi-\gp\}$ for some $0<\gp<\pi/2$. 
Our main aim is to describe the pseudodifferential structure of the resolvent
$(A-\gl)^{-1}$ as well as the asymptotic behavior, as $t\to 0^+$, of the  
operator $\etA$ and its trace $\tr\etA$ (heat trace). 

{}From the analytic point of view a cone is a product $(0,c)\times X$ 
together with a metric of the form $dr^2+r^2 g_X(r)$, where $g_X(r)$ is a 
smooth family of Riemannian metrics on the `cone base' $X$. Here, $X$ is
assumed to be a smooth compact manifold without boundary. 
For this reason, the analysis on a manifold with conical singularities 
takes place on a manifold with boundary $\mathbb{B}$ with the mentioned 
product structure near $\partial\mathbb{B}=X$. 
The natural differential operators appearing in this context are given, near 
the boundary, by the so-called operators of Fuchs type, cf. equation
\eqref{diff-op}. These operators arise, for instance, when  
considering Laplace-Beltrami operators associated to the metrics given 
above. A {\em cone differential operator} is a differential operator on the 
manifold $\mathbb{B}$ which is of Fuchs type near the boundary
$\partial\mathbb{B}$. 

There is a large number of papers concerning asymptotic expansions for 
resolvents, heat kernels and heat traces on smooth compact manifolds with 
and without boundaries, and their applications to geometry, index theory, 
mathematical physics and other areas. A large collection of references 
related to this topic can be found in {\sc Gilkey}'s book
\cite[Chapter~5]{Gilkey}. Let us especially mention the important papers 
\cite{ABP73}, \cite{Gilk73}, \cite{Gre71} and \cite{MP}.
An elegant pseudodifferential method for the study of resolvents of 
elliptic differential operators on smooth compact manifolds was introduced by
{\sc Seeley} in \cite{See67}, \cite{See69a}, \cite{See69b}. The idea of his
approach is to approximate the resolvents by means of parametrices that are
constructed  within a corresponding parameter-dependent pseudodifferential
calculus, where the parameter (up to a natural anisotropy) is involved as an
additional covariable. For genuine pseudodifferential operators the situation
is more delicate and other methods are required, cf. \cite{DG75},
\cite{Grubb}, \cite{GrSe95}, \cite{W1}.

On manifolds with conical singularities the analysis becomes more complicated.
However, for certain classes of self-adjoint operators on some singular
manifolds there are many results concerning resolvent, heat kernel and heat
trace asymptotics; in this context we want to refer to \cite{BrSe85}, 
\cite{BrSe89}, \cite{Ca83}, \cite{Ch79}, \cite{Ch83}, \cite{F1}, \cite{Kar98},
\cite{Le97} and \cite{Moo96}. In particular, {\sc Lesch} \cite{Le97} 
generalized the method introduced by {\sc Br\"uning} and {\sc Seeley} in
\cite{BrSe85} and obtained asymptotic expansions of the heat trace for
self-adjoint differential operators with coefficients that are independent of
the radial variable near the singularities. For the more general dependent 
case he obtained partial expansions. Let us finally mention that 
{\sc Loya} \cite{LoRes01}, \cite{LoHeat01} has recently obtained similar 
results to those from this paper by using geometric
blow-up techniques of {\sc Melrose}.
 
The strategy of this work is to follow the resolvent approach mentioned above
by introducing a suitable notion of parameter-dependent ellipticity 
(generalizing self-adjointness), and by giving an explicit construction 
of a parametrix within a natural class of operator-valued symbols.
We also take advantage of a parameter-dependent symbol class introduced by 
{\sc Grubb} and {\sc Seeley} \cite{GrSe95} in order to consider successfully
operators whose coefficients do depend smoothly on the radial variable $r$ 
near $0$ (i.e., near $\partial\mathbb{B}$), cf. equation \eqref{diff-op}.
  
In Section~\ref{coneBasics} we introduce the class of cone differential 
operators and discuss briefly some of the typical elements of 
{\sc Schulze}'s cone calculus. For an extensive study of cone algebras we 
refer to \cite{Sz91}, \cite{Sz98} and \cite{ScSzII}. 
In the third section we analyze the operator family $A-\gl$ for a cone 
differential operator $A$. The main idea is to consider $A-\gl$ as an element
of a certain class of operator-valued symbols, defined in an abstract setting 
by means of strongly continuous groups acting on Banach spaces. The parameter
$\gl$ plays the role of the covariable and is treated anisotropically with
respect to the covariables of the local symbols. Moreover,
by freezing the coefficients of $A$ at the boundary, cf. equation 
\eqref{pdell3}, we get a canonical object (the twisted homogeneous principal
symbol) that together with the symbolic structure of $A$ characterizes the
parameter-ellipticity of $A-\gl$. Our concept of ellipticity turns out to be
sufficient and necessary for the invertibility of $A-\gl$ on the canonical 
weighted Sobolev spaces. Next, we construct a {\em parametrix}, that is,
an element of the class that inverts $A-\gl$ modulo operator-valued
functions decreasing rapidly in the parameter and taking values in the Green
cone operators. This parametrix construction essentially relies on techniques
introduced by {\sc Schulze} \cite{Sz89} for the analysis on manifolds with
edges. In Section~\ref{weak-exp} we consider a modified version of the
parameter-dependent pseudodifferential calculus developed in \cite{GrSe95}.
This class contains the holomorphic symbols arising from the Mellin
pseudodifferential theory considered in this work.
The main purpose is to achieve an asymptotic expansion of the
local symbols in the parameter $\gl$ as $|\gl|\to\infty$. 

The last section is devoted to the study of the heat operator associated to
a cone differential operator $A$ of order $m>0$ which is parameter-elliptic 
in a suitable sector $\gL\subset\C$. Using the resolvent 
estimates obtained here, we define the semigroup $\{\etA\}_{t>0}$
by means of Dunford integrals. It is then proved that for each $t>0$, 
the operator $\etA$ belongs to the class of Green elements in the cone 
algebra (see Appendix~\ref{s-cone}). Green operators in the cone theory are
actually integral operators whose kernels are smooth in the interior
and behave in a particular way near the boundary. As a consequence, each
operator $\etA$ ($t>0$) is of trace class on the weighted Sobolev space
$\mathcal{H}^{s,\gg}(\mathbb{B})$ for all $s\in\R$, whenever $A$ is
parameter-elliptic with respect to the weight $\gg$.
Moreover, we show that the heat trace admits the expansion 
\[ \tr e^{-tA}\sim \sum_{k=0}^\infty C_k\,t^{(k-n-1)/m}+
   \sum_{k=0}^\infty C_k'\,t^{k/m}\log{t},\quad\text{as }\, t\to 0^+, \]
where $n=\dim\partial\mathbb{B}$. The coefficients $C_k$ and $C_k'$ depend on
the principal symbol of the cone operator $A$, and on its boundary spectrum
(cf. Section~\ref{cone-do}).

Roughly speaking, the construction relies on the use of three different
pseudodifferential calculi of local and global nature. 
The problem is reduced to an asymptotic expansion of the resolvent power
$(A-\gl)^{-\ell}$, for $\ell\in\N$ sufficiently large, which is viewed as
an operator-valued symbol of negative order. Using the standard symbolic
calculus we get a first approximation of $(A-\gl)^{-\ell}$ by means of 
a parametrix of $(A-\gl)^\ell$. The part of this parametrix localized in the
interior of $\mathbb{B}$ provides a `good' approximation of the resolvent 
since the local symbols are homogeneous. In order to achieve a suitable 
approximation of the resolvent near the boundary, we refine the parametrix 
by decomposing it into twisted homogeneous operators, separating at the same 
time the smoothing Green terms from the nonsmoothing parameter-dependent 
Mellin operators with degenerate symbols. The twisted homogeneity, which 
induces scalability of the kernels, permits us to imitate in a global sense 
the usual homogeneity arguments, but it is not enough to obtain the expansion
of the nonsmoothing Mellin components. Nevertheless, a complete asymptotic 
expansion of these operators can be achieved making use of the weakly 
parametric calculus from Section~\ref{weak-exp}. 

Finally, in order to reduce the search through the literature, 
some basic results and definitions are summarized in the Appendix.

\section{Elements of  the cone algebra}\label{coneBasics}
\subsection{Cone differential operators}\label{cone-do}
Let $X$ be a smooth manifold of dimension $n$. Denote by $X^\wedge$
the space $\R_+\times X$. A differential operator in $\diff^{m}(X^\wedge)$ 
is said  to be of {\em Fuchs type} if, expressed in the coordinates
$(r,x)$, it is of the form   
\begin{equation}\label{diff-op}     	
   A=r^{-m}\sumk{m} a_k(r)(-r\partial_r)^k   
\end{equation} 
with $a_k\in C^{\infty}(\rpbar,\diff^{m-k}(X))$. $\diff^\nu$ denotes the
space of differential operators of order $\nu\in\N_0$ with
smooth coefficients.

\begin{example}\label{fuchs2}
Let $X^\wedge$ be equipped with the {\em cone metric}
$dr^2+r^2g_{_X}(r)$, where $g_{_X}(r)$ is a family of Riemannian metrics on
$X$, smooth in $r\in\rpbar=[0,\infty)$.\\
The Laplace-Beltrami operator $\Delta_{\mathrm{cm}}$ corresponding 
to this metric is an operator of Fuchs type. In fact,  
\[ \Delta_{\mathrm{cm}} =r^{-2}\left\{\Delta_{g_{_X}\!(r)}
   +\left\{-n+1-rG^{-1}(\partial_r G)\right\}(-r\partial_r)
   +(-r\partial_r)^2\right\}, \]
where, in local coordinates $(r,x_1,\dots,x_n)$, 
$G(r,x)=|\det(g_{_X}(r)(\partial_{x_i},\partial_{x_j}))|^{1/2}$.
\end{example}

On a smooth manifold with boundary $M$ we may consider the so-called
{\em $b$-tangent bundle} ${}^b TM$ (cf. {\sc Melrose}
\cite[Section~2.2]{Mel93}). For a differential operator of Fuchs type $A$
as in \eqref{diff-op} there is a function $\gs_{\psi,b}^m(A)$ in
$C^{\infty}({}^b T^*X^\wedge\onull)$ such that, in local coordinates
$(r,x)\in\R_+\times U$, $U\subset X$, it takes the form
\begin{equation}\label{fuchsS}
\gs_{\psi,b}^m(A)(r,x,\gr,\xi)= r^m \gs_{\psi}^m(A)(r,x,\gr/r,\xi),
\end{equation} 
where $\gs_{\psi}^m(A)$ is the usual homogeneous principal  symbol of $A$
on $X^\wedge$. ${}^b T^*X^\wedge$ denotes the dual of ${}^b TX^\wedge$, 
and $(\gr,\xi)\in\R\times\R^n$ are the covariables to $(r,x)$.

Fredholm properties of Fuchs type operators are determined by a pair of
symbols. The first one is the homogeneous principal symbol 
$\gs_{\psi,b}^m(A)$ which characterizes the ellipticity in the interior. 
The second one is the so-called conormal symbol of $A$ which is an
operator-valued symbol, living at the boundary, that can be described by 
means of the {\em Mellin transform}
\[ \mathcal{M}_{r\to z}u(z)=\int_0^\infty r^{z-1} u(r) dr 
   \;\text{ for } u\in C^{\infty}_{0}(\R_+). \] 
The point is that the totally characteristic derivative $(-r\partial_r)$
corresponds, in the image of the Mellin transform, to multiplication by the
complex variable $z$. Therefore, if $a_k\in C^{\infty}(\rpbar,\diff^{m-k}(X))$
for $k=0,\dots,m$, are the coefficients of the operator $A$ from
\eqref{diff-op}, then the operator-valued polynomial
\begin{equation}\label{f2m.1}  
 h(r,z)=\sumk{m}a_k(r)z^k\quad \text{with }\,(r,z)\in\rpbar\times\C,  
\end{equation} 
may be interpreted as the corresponding {\em Mellin symbol} of $A$. 
In other words, for any real $\gb$ and $u\in C^{\infty}_{0}(X^\wedge)=
C^{\infty}_{0}(\R_+,C^{\infty}(X))$ we have the relation 
\begin{equation}\label{f2m.2} 
  Au(r)=[\op_M(h)u](r):=\frac{1}{2{\pi}i}
 \int_{\gG_{\gb}}r^{-z-m}h(r,z)[\mathcal{M}_{r'\to z}u](z)dz,  
\end{equation}
where $\gG_\gb:=\{z\in\C\sod \Re z=\gb\}$. Observe that $z\mapsto 
h(r,z):\C\to\diff^m(X)$ is holomorphic for each $r\in\R_+$. 
Now, the {\em conormal symbol} $\sM^m(A)$ of the Fuchs type operator $A$
is just the polynomial \eqref{f2m.1} evaluated at $r=0$, i.e.,
\[ \sM^m(A)(z):=h(0,z)=\sumk{m}a_k(0)z^k. \] 
The property of being Fuchs type is invariant under changes of
coordinates, so \eqref{diff-op} serves as model for a class of operators on
manifolds with boundaries. After a suitable blow-up, a manifold with conical
singularities becomes a smooth manifold with boundary with a structure that 
can be precisely described by operators of Fuchs type.

\begin{definition}\label{fcone1} 
Let $\mathbb{B}$ be a smooth compact manifold with boundary
$\partial\mathbb{B}=X$. The space $\diffC{m}$ of {\em cone differential
operators} consists of all operators in $\diff^m(\iBB)$ which are of 
Fuchs type near the boundary. 
\end{definition}

Fix once and for all a defining function $\g$ for the boundary of $\mathbb{B}$,
that is, a smooth function $\g:\mathbb{B}\to\rpbar$ such that $\g$ is positive
in the interior of $\mathbb{B}$, it vanishes on $\partial\mathbb{B}$, and
$d\g\not=0$ on $\partial\mathbb{B}$. We assume that $\g^{-1}([0,2))$ is a
collar neighborhood of $\partial\mathbb{B}$ in $\mathbb{B}$.

The conormal symbol of a cone operator $A$ is defined as the
operator family
\[ u\mapsto \g^{m-z} A(\g^z\tilde u)|_{\g=0}:C^{\infty}(X)\to C^{\infty}(X),\]
where $\tilde u$ is some extension of $u$.

A function $\go\in C^{\infty}_{0}(\rpbar)$ is called a (boundary) 
{\em cut-off function} if $\supp\go\subset[0,2)$ and $\go=1$ near $r=0$. 
Note that every $\go$ can be viewed as a function on $[0,2)\times X$ as 
well as a function on $\mathbb{B}$, extending by zero. We say that two
functions $\phi$, $\psi$ satisfy the relation $\phi\prec\psi$ if and only 
if $\phi\psi=\phi$.

\begin{remark}
Let $\go_0$, $\go_1$ and $\go_2$ be cut-off functions with 
$\go_2\prec\go_1\prec\go_0$. Let $\widetilde{\mathbb{B}}$ be the double of
$\mathbb{B}$. Since differential operators are local, i.e., $\supp
Au\subset\supp u$, the elements of $\diffC{m}$ can be written in the form 
\[ A=\go_1 A_0\,\go_0+(1-\go_1) A_1(1-\go_2), \] 
where $A_0$ is of Fuchs type on $X^\wedge$ and
$A_1\in\diff^m(\widetilde{\mathbb{B}})$.
In this decomposition, $A_0$ is not canonical but its conormal symbols 
is that of the operator $A$. Further, the interior principal symbol 
$\gs_{\psi}^m(A)$ of $A$ induces a symbol 
$\gs_{\psi,b}^m(A)\in C^{\infty}({}^b T^*\mathbb{B}\onull)$ that, in
local coordinates near the boundary, satisfies the relation \eqref{fuchsS}.
\end{remark}
\begin{proposition}\label{fcone2}  
The class of cone differential operators is closed under compositions. 
Moreover, the conormal symbols satisfy the relation
\begin{equation}\label{comp-co}   
 \sM^{m_1+m_2}(A_2 A_1)(z)=\sM^{m_2}(A_2)(z+m_1)\sM^{m_1}(A_1)(z). 
\end{equation}
\end{proposition}
\begin{definition}\label{fcone3}  
$A\in\diffC{m}$ is called {\em elliptic with respect to $\gg\in\R$} if 
\begin{enumerate} 
\item[(i)] $\gs_{\psi,b}^m(A)\not=0$ on ${}^b T^*\mathbb{B}\onull$,
\item[(ii)] $\sM^m(A)(z):H^m(X)\to L^2(X)$ is an isomorphism for all
$z\in\gG_{\frac{n+1}{2}-\gg}$.
\end{enumerate}
\end{definition}

As in \cite{MM} we denote
\[ \spec_b(A):=\{z\in\C\sod \sM^m(A)(z):H^m(X)\to L^2(X) 
   \text{ is \underline{not} an isomorphism}\}. \] 
This set is called the {\em boundary spectrum} of $A$ and is known to be
discrete, and finite on vertical strips. 

\begin{example}\label{fcone4}
Let $g$ be a metric on $\mathbb{B}$ that coincides with a cone metric
in a collar neighborhood of $\partial\mathbb{B}=X$, cf. Example~\ref{fuchs2}. 
Then the corresponding  Laplace-Beltrami operator $\Delta_g$ is an 
elliptic cone differential operator with conormal symbol
\[ \sM^2(\Delta_g)(z)=\Delta_{g_{_X}\!(0)}-(n-1)z+z^2,\quad n=\dim X. \]
In this case, $\, z\in\spec_b(\Delta_g)$ if and only if 
$\,(n-1)z-z^2\in\spec_{L^2(X)}(\Delta_{g_{_X}\!(0)})$. 
\end{example}

\subsection{Weighted Sobolev spaces} \label{sobolev}
Let us introduce a scale of weighted (cone) Sobolev spaces on which the cone
operators act continuously. These spaces are defined in a similar way as the
usual Sobolev spaces but based on the Mellin (instead of the Fourier)
transform in the singular direction. 
Recall that $\gG_\gb=\{z\in\C\sod \Re z=\gb\}$.

\smallskip
For $s,\gg\in\R$ let $\hsg(\R_+\times\R^n)$ be the closure of 
$C^{\infty}_{0}(\R_+\times\R^n)$ with respect to
\[ \norm{u}_{\hsg}^2 =\frac{1}{2\pi i}
   \iint_{\gG_{\frac{n+1}{2}-\gg}}\!\!(1+|z|^2+|\xi|^2)^s\,
   |(\mathcal{M}_{r\to z}\mathcal{F}_{x\to\xi}u)(z,\xi)|^2\,dz\,d\xi. \] 
The transformation $\,S_{\gg-\frac{n}{2}}:
C^{\infty}_{0}(\R_+\times\R^n)\to C^{\infty}_{0}(\R^{1+n})$ given by 
\[ \left(S_{\gg-\frac{n}{2}}\,u\right)(r,x) =
   e^{-(\frac{n+1}{2}-\gg)r}u(e^{-r},x) \]
extends to an isomorphism  $\,S_{\gg-\frac{n}{2}}:\hsg(\R_+\times\R^n)
\to H^s(\R^{1+n})$, where $H^s$ denotes the standard Sobolev space. 
We get the equivalence of norms 
\begin{equation}\label{MF-equiv}  
  \norm{u}_{\hsg(\R_+\times\R^n)}
  \sim \big\|{S_{\gg-\frac{n}{2}}u}\big\|_{H^s(\R^{1+n})}.
\end{equation}

In order to define the suitable Sobolev spaces on the manifold $X^\wedge$ 
we fix an open covering $\{U_1,\ldots,U_N\}$ of $X$ with corresponding 
diffeomorphisms $\chi_j:U_j\to V_j\subset\R^n$ and 
$\tilde\gk_j:U_j\to\tilde V_j\subset S^n$, where $S^n$ is the unit sphere in 
$\R^{1+n}$. In addition, we define $\gk_j:\R_+\times 
U_j\to\R^{1+n} \setminus\{0\}$ by $\gk_j(r,x):=r\tilde\gk_j(x)$.
Let finally $\{\phi_1,\ldots,\phi_N\}$ be a subordinate partition of unity.

\begin{definition}\label{sob2}
Let $\go$ be a cut-off function. For $s,\gg\in\R$ let 
$\mathcal{K}^{s,\gg}(X^\wedge)$ denote the closure of
$C^{\infty}_{0}(X^\wedge)$ with respect to the norm
\begin{equation}\label{knorm} 
 \norm{u}_{\mathcal{K}^{s,\gg}}^2= \sum_{j=1}^N \Big\{ 
 \|\chi_{j*}(\phi_j\,\go u)\|_{\hsg(\R_+\times\R^n)}^2 +
 \|\gk_{j*}(\phi_j(1-\go)u)\|_{H^{s}(\R^{1+n})}^2 \Big\}.
\end{equation} 
\end{definition}

Note that for every fixed localization data (partition of unity, 
diffeomorphisms $\chi_j$ and $\tilde\gk_j$, and cut-off function $\go$) the
space $\mathcal{K}^{s,\gg}(X^\wedge)$ has a Hilbert space structure. 
Another choice of data yields a norm equivalent to \eqref{knorm}.  
The spaces $\mathcal{K}^{s,\gg}(X^\wedge)$ are subspaces of
$H_{loc}^s(X^\wedge)$. In particular, we have $\mathcal{K}^{0,0}(X^\wedge)
= L^2(X^\wedge, r^n dr dx)$.

\begin{proposition}\label{fcont1} 
Let $A\in\diff^{m}(X^\wedge)$ be of Fuchs type with coefficients 
$a_k(r)$ that are independent of $r$ for large values of $r$. 
Then the map
\[ A:\mathcal{K}^{s,\gg}(X^\wedge)\to \mathcal{K}^{s-m,\gg-m}(X^\wedge)
   \text{ is continuous for all } s,\gg\in\R. \]
\end{proposition}

\begin{definition}\label{sob3} 
For $s,\gg\in\R$ and $\go$ a cut-off function on $\mathbb{B}$,  let 
\[ \hsg(\mathbb{B}):=\{u\in\mathcal{D}'(\iBB)\sod
   \go u\in\mathcal{K}^{s,\gg}(X^\wedge)\text{ and } 
   (1-\go)u\in H^s(\widetilde{\mathbb{B}})\}.  \]
\end{definition}

We endow this space with the norm 
\[ \norm{u}_{\hsg(\mathbb{B})}=\norm{\go u}_{\mathcal{K}^{s,\gamma}(X^\wedge)}
  + \norm{(1-\go)u}_{H^s(\widetilde{\mathbb{B}})}. \]
Another cut-off function $\go$ leads to an equivalent norm. 
$\hsg(\mathbb{B})$ is a subspace of $H_{loc}^s(\iBB)$ and
$\mathcal{H}^{0,0}(\mathbb{B})\cong\g^{-n/2}L^2(\mathbb{B})$ for any boundary
defining function $\g$. Further, $\mathcal{H}^{s,\gg}(\mathbb{B}) =
\g^\gg\mathcal{H}^{s,0}(\mathbb{B})$ for every $\gg\in\R$. 
Moreover, the space $\hsg(\mathbb{B})$ is a Hilbert space.

\begin{lemma}\label{op-norm}
Let $s,s',\gg,\gg'\in\R$, and let $\go_0$, $\go_1$ and $\go_2$ be cut-off
functions. Let $A_0$ and $A_1$ be operators on $X^\wedge$ and
$\widetilde{\mathbb{B}}$, respectively, such that 
\[ \go_1 A_0\,\go_0 \in\mathcal{L}(\mathcal{K}^{s,\gg},\mathcal{K}^{s',\gg'})
  \;\text{ and }\;(1-\go_1) A_1(1-\go_2) \in\mathcal{L}(H^{s},H^{s'}). \]
If $A=\go_1 A_0\,\go_0+(1-\go_1) A_1(1-\go_2)$, 
then there exists $C>0$ such that
\begin{multline*}
 \norm{A}_{\mathcal{L}(\mathcal{H}^{s,\gg}(\mathbb{B}),
 \mathcal{H}^{s',\gg'}(\mathbb{B}))}\\
 \leq C\left(\norm{\go_1 A_0\,\go_0}_{\mathcal{L}(\mathcal{K}^{s,\gg},
 \mathcal{K}^{s',\gg'})} + \norm{(1-\go_1)
 A_1(1-\go_2)}_{\mathcal{L}(H^{s},H^{s'})}\right). 
\end{multline*}
\end{lemma}

\begin{proposition}\label{fcont2} 
Any $A\in\diffC{m}$ extends to a continuous operator
\[ A:\mathcal{H}^{s,\gg}(\mathbb{B})\to \mathcal{H}^{s-m,\gg-m}(\mathbb{B}) 
   \;\text{ for all } s,\gg\in\R.\]
\end{proposition}

{\bf Fredholm property.}
The concept of ellipticity introduced in  Definition~\ref{fcone3} leads to
the Fredholm property of the elliptic operators on the weighted Sobolev spaces
on $\mathbb{B}$. More precisely, if  $A\in\diffC{m}$ is elliptic with respect
to $\gg$, then the map
\[ A:\hsg(\mathbb{B})\to\mathcal{H}^{s-m,\gg-m}(\mathbb{B})
  \text{ is Fredholm for all } s\in\R. \]
This is a consequence of the following 

\begin{theorem}\label{cone5}
To every $A\in\diffC{m}$, elliptic with respect to $\gg$, there exists a 
cone pseudodifferential operator $Q$ of order $-m$ such that
\[ QA-1:\mathcal{H}^{s,\gg}(\mathbb{B})\to 
   \mathcal{H}^{\infty,\gg}(\mathbb{B}) 
   \;\text{ and }\; AQ-1:\mathcal{H}^{s,\gg-m}(\mathbb{B})\to
   \mathcal{H}^{\infty,\gg-m}(\mathbb{B}) \]
are smoothing operators of Green type (Appendix~\ref{s-cone}). 
$Q$ is called a {\em parametrix} of $A$.
\end{theorem}

\begin{remark} 
If $A$ is elliptic with respect to two different weights, say $\gg_1$ and
$\gg_2$, then the difference between the Fredholm indices depends on the
elements in the boundary spectrum of $A$ situated between the
corresponding weight lines in $\C$. Moreover, there is an explicit formula 
for the difference $\text{ind}_{\gg_1}A -\text{ind}_{\gg_2}A$,
cf. \cite{MM}, \cite{Sz91}. 
\end{remark}

{\bf Adjoint operators.}
The restriction to $C^{\infty}_{0}(X^\wedge)\times C^{\infty}_{0}(X^\wedge)$ 
of the sesquilinear pairing $(\cdot,\cdot)_{\mathcal{K}^{0,0}}$ extends to a
nondegenerate pairing 
\[ \mathcal{K}^{s,\gg}(X^\wedge)\times \mathcal{K}^{-s,-\gg}(X^\wedge)\to\C 
  \;\text{ for any } s,\gg\in\R. \] 
The {\em formal adjoint} of $A\in\mathcal{L}(\mathcal{K}^{s,\gg},
\mathcal{K}^{s-m,\gg-m})$ with respect to $(\cdot,\cdot)_{\mathcal{K}^{0,0}}$ 
is the unique operator $A^*_{0}$ satisfying $(Au,v)=(u,A^*_{0}v)$ for all
$u,v\in C^{\infty}_{0}(X^\wedge)$. More precisely, if $A$ is a Fuchs type
operator as in \eqref{diff-op}, and $u\in C^{\infty}_{0}(X^\wedge)$, then
\[ A^*_0 u(r)=r^{-m}\sumk{m}(n+1-m+r\partial_r)^k\{a_k(r)^\star u(r)\}, \] 
where $a_k(r)^\star$ denotes the pointwise formal adjoint in $L^2(X)$. 
For an arbitrary $\gb\in\R$ the adjoint of $A$ with respect to
$(\cdot,\cdot)_{\mathcal{K}^{0,\gb}}$ is given by
$A^*_\gb=\g^{2\gb}A^*_0\,\g^{-2\gb}$. Here $\g$ is a defining function for 
the boundary of $X^\wedge$ with $\g=1$ near infinity.

\section{Parameter-dependent calculus}\label{p-cone}
In this section we basically concentrate our attention to the operator
family $A-\gl$ as the typical parameter-dependent cone operator.
In fact, all the results exposed here can be generalized to wider classes
of parameter-dependent cone pseudodifferential operators, cf. \cite{BuSz97},
\cite{Gil}, \cite{ScSzES}, \cite{Sz98}, \cite{Seiler}.

\smallskip
Let $\gL\subset\C$ be a closed angle with vertex at the origin.
Let $X$ be a smooth compact manifold, $\dim X=n$, and let $\mathbb{B}$ be 
the manifold with boundary $\partial\mathbb{B}=X$ introduced in
Section~\ref{cone-do}. Recall that $\widetilde{\mathbb{B}}$ denotes the 
double of $\mathbb{B}$. Again, let $X^\wedge=\R_+\times X$. 

\subsection{Ellipticity condition for $\mathbf{A-\lambda}$}
If $A$ is a cone differential operator of order $m\in\N$, and
$\gl\in\gL$, then $A-\gl$ can be written in the form
\begin{equation}\label{pdell1} 
 A-\gl=\go_1\Big(r^{-m}\sumk{m}a_k(r)(-r\partial_r)^k-\gl\Big)\go_0      
 +(1-\go_1)(P-\gl)(1-\go_2),  
\end{equation}
where $a_k\in C^{\infty}(\rpbar,\diff^{m-k}(X))$,
$\,P\in\diff^m(\widetilde{\mathbb{B}})$, and $\go_0$, $\go_1$ and $\go_2$ are
cut-off functions with $\go_2\prec\go_1\prec\go_0$. Denote $a(\gl)=A-\gl$ and 
\begin{equation}\label{pdell2}
 a_{\wedge}(\gl)= r^{-m}\sumk{m} a_k(r)(-r\partial_r)^k-\gl
 \in L_{\cl}^{m,m}(X^\wedge;\gL).
\end{equation}
Because of the presence of $\go_0$ and $\go_1$ in \eqref{pdell1}
we may assume that every $a_k(r)$ is independent of $r$ for $r>2$. 
Then, for every $s,\gg\in\R$ we have 
\[ a_{\wedge}\in S^{m,m}(\gL;\mathcal{K}^{s,\gg}(X^\wedge),
   \mathcal{K}^{s-m,\gg-m}(X^\wedge)) 
   \quad(\text{cf. Appendix~\ref{opv-symb}}) \]
with the strongly continuous group of isomorphisms defined by
\[ (\gk_{\tau} u)(r,x):={\tau}^{\frac{n+1}{2}} u({\tau}r,x)\;
   \text{ for }\tau>0 \,\text{ and } n=\dim X. \]
A proof of this statement can be found in \cite[Prop. 3.1.1]{Sz99} or 
\cite[Prop. 3.3.38]{Sz98}. \\
Finally, let 
\begin{equation}\label{pdell3}
\sK^m(a)(\gl):= r^{-m}\sumk{m} a_k(0)(-r\partial_r)^k-\gl.
\end{equation}
This operator-valued symbol satisfies the homogeneity relation
\[ \sK^m(a)(\tau^m\gl)=\tau^m\gk_{\tau}\sK^m(a)(\gl)\gk_{\tau}^{-1} \] 
for all $\tau>0$, $\gl\in\gL\setminus\{0\}$.
For this reason, \eqref{pdell3} will be called the {\em twisted homogeneous
principal symbol} of $a(\gl)$. Due to the local property of the operator $A$,
$\sK^m(a)$ is well-defined and certainly a canonical object that can be used 
to define the parameter-dependent ellipticity of $A-\gl$:

\begin{definition}\label{eell1}  
The family $a(\gl)=A-\gl$ is said to be {\em parameter-elliptic}
with respect to $\gg\in\R$ on $\gL$ if and only if 
\begin{enumerate}
\item[(i)]
 $\gs_{\psi,b}^m(A)-\gl\not=0\,$ on $\,({}^bT^*\mathbb{B}\times\gL)\onull$, 
\item[(ii)]
 $\spec_b(A)\cap\gG_{\frac{n+1}{2}-\gg}=\emptyset$,
 \quad (cf. Definition~\ref{fcone3})
\item[(iii)]
 $\sK^m(a)(\gl):\mathcal{K}^{s,\gg}(X^\wedge) \to 
 \mathcal{K}^{s-m,\gg-m}(X^\wedge)$ is an isomorphism for some $s\in\R$ 
 and every $\gl\in\gL$ sufficiently large.
\end{enumerate}
\end{definition}

\begin{theorem}\label{einv1} 
If $A-\gl$ is parameter-elliptic with respect to $\gg$, then there is a
constant  $R>0$ such that $A-\gl:\mathcal{H}^{s,\gg}(\mathbb{B})\to
\mathcal{H}^{s-m,\gg-m}(\mathbb{B})$ is invertible 
for every $\gl\in\gL$ with $|\gl|\ge R$, and all $s\in\R$. Furthermore,
there are constants $C(s,\gg,m)>0$ and $M(s,\gg,m)\ge 0$ such that
\[ \norm{(A-\gl)^{-1}}_{\mathcal{L}(\mathcal{H}^{s-m,\gg-m})}
   \leq C\,(1+|\gl|)^{-1+M/m}. \]
\end{theorem}
\begin{proof} 
In Theorem~\ref{eparam} we will construct a parametrix $b(\gl)$
of $a(\gl)=A-\gl$ such that $g=ab-1\in\mathcal{S}(\gL,\coneg{\gg-m,\gg-m})$.
Moreover, this construction will provide a family $\{b(\gl)\}_{\gl\in\gL}
\subset \mathcal{L}(\mathcal{H}^{s-m,\gg-m}(\mathbb{B}),
\mathcal{H}^{s,\gg}(\mathbb{B}))$ that belongs to the class
\[ S_{\cl}^{-m+M,m}(\gL;\mathcal{H}^{s-m,\gg-m}(\mathbb{B}),
  \mathcal{H}^{s-m,\gg-m}(\mathbb{B})) \] 
for some $M=M_{s,\gg,m}\ge 0$, cf. Remark~\ref{opv-symbol} and 
Appendix~\ref{opv-symb}.

Fix now $s_0\in\R$. In $\mathcal{L}(\mathcal{H}^{s_0-m,\gg-m}(\mathbb{B}))$
we have the estimate 
\[ \norm{a(\gl)b(\gl)-1}\leq C(1+|\gl|)^{-N}\;\text{ for every } N\in\N.\]
Hence for some $R>0$ we get $\norm{a(\gl)b(\gl)-1}\leq 1/2$ 
for all $|\gl|\geq R$. Therefore, $a(\gl)b(\gl)$ is invertible in 
$\mathcal{L}(\mathcal{H}^{s_0-m,\gg-m}(\mathbb{B}))$ for  $|\gl|\geq R$, and 
$\norm{[a(\gl)b(\gl)]^{-1}}\leq 2$. In fact, the inverse of $a(\gl)$ is 
given by $b(\gl)[a(\gl)b(\gl)]^{-1}$. Furthermore, 
\[ \norm{a(\gl)^{-1}}=\norm{b(\gl)[a(\gl)b(\gl)]^{-1}}
  \leq 2\norm{b(\gl)}\leq C_{s_0,\gg,m}\,(1+|\gl|)^{-1+M/m} \]
for some constant $C_{s_0,\gg,m}$. The last inequality is the symbol estimate
of $b(\gl)$ as an operator-valued symbol of order $(-m+M,m)$.
In order to see that for $|\gl|\geq R$ the operator $a(\gl)$ is invertible 
even for all $s\in\R$, let $\tilde g=ba-1\in\mathcal{S}(\gL,\coneg{\gg,\gg})$.
Then the inverse of $a(\lambda)$ can also be written as
\[ a(\gl)^{-1}=b(\gl)-b(\gl)g(\gl)+{\tilde g(\gl)}a(\gl)^{-1}g(\gl) \]
which belongs to $\mathcal{L}(\mathcal{H}^{s-m,\gg-m}(\mathbb{B}),
\mathcal{H}^{s,\gg}(\mathbb{B}))$ for all $s\in\R$ and $|\gl|\geq R$. 
Observe that the inverse of $a(\gl)$ appearing between $\tilde g$ and $g$ 
is the inverse with respect to $s_0$.
\end{proof}

\begin{corollary}\label{einv2}
If $A-\gl$ is parameter-elliptic with respect to $\gg$ and $b(\gl)$ 
is a corresponding parametrix, then there exists $R>0$ such that
\[ (A-\gl)^{-1}-b(\gl)\in\mathcal{S}(\gL_R,\coneg{\gg-m,\gg}),\]
where $\gL_R=\{\gl\in\gL \sod |\gl|\ge R\}$.  
\end{corollary}

Although the parameter-ellipticity introduced above seems to be a  
strong condition, it is indeed necessary for the invertibility of 
$A-\gl$ on the weighted Sobolev spaces.

\begin{theorem}\label{einv3}
Let $A\in\diffC{m}$ be such that $A-\gl:\mathcal{H}^{s,\gg}(\mathbb{B})\to
\mathcal{H}^{s-m,\gg-m}(\mathbb{B})$ is invertible for all $\gl\in\gL$
sufficiently large, and $\norm{(A-\gl)^{-1}}$ is uniformly bounded. 
Then $A-\gl$ is parameter-elliptic with respect to $\gg$ in the sense of 
Definition~\ref{eell1}.
\end{theorem}

A proof of this theorem together with other results concerning resolvents 
of cone operators on arbitrary domains will be given in \cite{GiMe1}.

\begin{example}\label{eell2}
Let $\gb\in C^{\infty}(X^\wedge)$ with $\gb(0)=0$. Let $D(r)$ be a smooth
family of elliptic selfadjoint differential operators on $X$ such that 
\[ \spec_{L^2(X)}(D(0))\subset (-\infty,\gd) \;\text{ with } 
 \gd=\big(\tfrac{n+1}{2}-|\gg-2|-2\big)\big(\tfrac{n+1}{2}-|\gg-2|\big).\]
Then the Fuchs type operator 
$A_0=r^{-2}\left\{D(r)-(1-n+\gb(r))(r\partial_r)+(r\partial_r)^2\right\}$
is a positive self-adjoint operator on $\mathcal{K}^{0,\gg-2}(X^\wedge)$ 
with $\mathcal{K}^{2,\gg}(X^\wedge)$ as domain. Further let
$A_1\in\diff^2(\widetilde{\mathbb{B}})$ be such that $A_1-\gl$ is
parameter-elliptic\footnote{In this case, parameter-elliptic on $\gL$ 
just means $\gs_{\psi}^m(A_1) -\gl\not=0$ on $T^*\widetilde{\mathbb{B}}
\times\gL\onull$.} on a conical set $\gL$ contained in the resolvent set 
of $A_0$. Then the operator family
\[ A-\gl:=\go_1(A_0-\gl)\go_0+(1-\go_1)(A_1-\gl)(1-\go_2), \quad\gl\in\gL,\]
is parameter-elliptic with respect to $\gg$ for any cut-off
functions $\go_2\prec\go_1\prec\go_0$.
\end{example}
\begin{example}\label{eell3}
Let $A\in\diffC{m}$ be such that $A-\gl$ is parameter-elliptic on 
$\gL$. If $B$ is a cone differential operator of order
less than $m$, then for any sufficiently small
$\eps$ the family $A+\eps B-\gl$ is also parameter-elliptic.
\end{example}

\subsection{Parametrix construction}\label{param-c}
Let us set again $a(\gl)=A-\gl$ for $\gl\in\gL$. Using \eqref{pdell1} and
\eqref{pdell2} we rewrite
\begin{equation}\label{pdell4}
  a(\gl)=\go_1\, a_{\wedge}(\gl)\go_0 +(1-\go_1)(P-\gl)(1-\go_2). 
\end{equation}
Our strategy to find a parametrix of $a(\gl)$ will be to construct a
parametrix of $a_{\wedge}(\gl)$ on $X^\wedge$ using techniques borrowed from
the edge symbolic calculus introduced by {\sc Schulze}, cf. \cite{Sz89},
\cite{Sz98}. The parametrix construction for $P-\gl$ is well-known and can be
found, for instance, in \cite[Section 11]{Shubin}.

\begin{definition}\label{edge2}  
Let $\mu\in\R$, $d\in\N$ and let $\bar d:=(1,d)$. 
Let $M_{\mathcal{O}}^{\mu,d}(X;\gL)$ be the space of holomorphic functions
\[ z\mapsto h(z,\gl)\in\mathcal{O}(\C,L^{\mu,d}(X;\gL))
    \quad\text{(cf. Appendix~\ref{pdo-p})} \] 
such that $h|_{\gG_\gb}\in L_{\cl}^{\mu,\bar d}(X;\gG_\gb\times\gL)$ 
for every $\gb\in\R$, uniformly for $\gb$ in compact intervals. 
\end{definition} 

The elements of $C^{\infty}(\rpbar,M_{\mathcal{O}}^{\mu,d}(X;\gL))$ will be
called {\em parameter-dependent holomorphic Mellin symbols}. The corresponding
smoothing class is defined by replacing $\mu$ by $-\infty$ and omitting $d$.
Recall that $\gG_\gb=\{z\in\C\sod \Re z=\gb\}$. 

\begin{example}\label{edge3}
The operator-valued function $h(r,z,\gl)=\sumk{m}a_k(r)z^k-r^m\gl$ is
clearly an element of $C^{\infty}(\rpbar,M_\mathcal{O}^{m,m}(X;\gL))$. 
\end{example}

Let $\coneg{\cdot,\cdot}$ denote the class of Green cone operators as defined
in Appendix~\ref{s-cone}, and let 
$\mathcal{S}\big(\gL,\coneg{\cdot,\cdot}\big)$ be the space of rapidly 
decreasing $\mathrm{C_G}$-valued functions. 
\begin{theorem} {\bf (Parametrix)} \label{eparam}
If $a(\gl)=A-\gl$ is parameter-elliptic with respect to $\gg$, 
there is a family $b(\gl)$ of cone pseudodifferential operators of order 
$-m$ such that 
\[ ba-1\in\mathcal{S}\big(\gL,\coneg{\gg,\gg}\big)\;\text{ and }\;
   ab-1\in\mathcal{S}\big(\gL,\coneg{\gg-m,\gg-m}\big).\]
The family $b(\gl)$ is a (parameter-dependent) parametrix of $a(\gl)$.  
\end{theorem}

To prove this theorem let us first construct a parametrix of $a_{\wedge}$ on 
$X^\wedge$, i.e., a pseudodifferential family $b_{\wedge}(\gl)$ such that  
\begin{equation}\label{pmtx1}
   b_\wedge a_{\wedge}-1\in\edgeg{-\infty}{\gg,\gg} \;\text{ and }\;
   a_{\wedge}b_{\wedge}-1\in\edgeg{-\infty}{\gg-m,\gg-m}, 
\end{equation}
where $\edgeg{-\infty}{\cdot,\cdot}$ is the class of Green operator-valued 
symbols from Appendix~\ref{opv-symb}.

The construction of $b_{\wedge}(\gl)$ will be carried out in several steps. 
\begin{lemma} {\bf (Step 1)}\label{param1}
There is a family $b_1(\gl)$ of pseudodifferential operators of order $-m$
such that
\begin{enumerate} 
\item \label{param1.1} 
 $a_{\wedge}(\gl)b_1(\gl)-1\in L^{-\infty}(X^\wedge;\gL)$,
\item \label{param1.2} 
 $\sM^0(a_{\wedge}b_1)(z)-1\in M^{-\infty}_{\mathcal{O}}(X)$,
\item \label{param1.3} 
 $a_{\wedge}b_1-1\in S_{\cl}^{0,m}(\gL;\mathcal{K}^{s,\gg-m}(X^\wedge),
 \mathcal{K}^{s',\gg-m}(X^\wedge)^\gs)$ for all $s,s',\gs\in\R$, where\\[1ex]
 $\mathcal{K}^{s',\gg-m}(X^\wedge)^\gs :=(1+r)^{-\gs}
 \,\mathcal{K}^{s',\gg-m}(X^\wedge)$.
\end{enumerate}
\end{lemma}
\begin{lemma} {\bf (Step 2)}\label{param2}
There is a meromorphic function $f$ such that if 
\[ v_1(\gl):=\go_1(r\ladm)\, r^{m}\op_M(f)\,\go_0(r\ladm) \]
with $[\gl]_m$ as in \eqref{opvs1.1},
then $a_{\wedge}(b_1+v_1)-1\in\edgeg{0,m}{\gg-m,\gg-m}$.
\end{lemma}
\begin{lemma} {\bf (Step 3)}\label{param3}
There is a Green symbol $g_1\in\edgeg{-m,m}{\gg-m,\gg}$ such that  
$b_{\wedge}^0:=b_1+v_1+g_1$ satisfies
$\sK^m(a_{\wedge})^{-1}=\sK^{-m}(b_{\wedge}^0)$. In particular, this implies
$a_{\wedge}b_{\wedge}^0 - 1\in\edgeg{-1,m}{\gg-m,\gg-m}$.
\end{lemma} 
\begin{lemma} {\bf (Step 4)}\label{param4}
There is a parametrix $b_{\wedge}$ of  $a_{\wedge}$ such that 
\eqref{pmtx1} holds.
\end{lemma} 

{\bf Proof of Lemma \ref{param1}.} 
In the representation \eqref{pdell1} we may choose the coefficients $a_k$ 
such that $a_k(r)=a_k(0)$ for $r\geq 2$. Thus $a_{\wedge}(\gl)$ may be assumed
to be parameter-dependent elliptic on $\gL$ as an operator family living on
$X^\wedge$. Therefore, there exists a family of pseudodifferential operators 
$b_0(\gl)\in L_{\cl}^{-m,m}(X^\wedge;\gL)$ such that 
\begin{equation}\label{roughparam}
 a_{\wedge}(\gl)b_0(\gl)-1\in L^{-\infty}(X^\wedge;\gL).  
\end{equation}
Moreover, $b_0(\gl)$ can be written as $r^m\op_{r,x}(q)(\gl)$ with a 
symbol $q$ such that
\[ q(r,x,\gr,\xi,\gl)=\tilde q(r,x,r\gr,\xi,r^m\gl)\;\text{ for some }
   \tilde q\in S_{\cl}^{-m,m}(\rpbar\times\gO\times\R^{1+n};\gL). \]
The so-called Mellin quantization theorem, see e.g. \cite[Theorem~3.2]{GSS2},
states that there is a Mellin symbol $h(r,z,\gl)= \tilde h(r,z,r^m\gl)$ with
$\tilde h\in C^\infty(\rpbar, M_{\mathcal{O}}^{-m,m}(X;\gL))$ such that 
$\op_{r,x}(q)(\gl)-\op_M(h)(\gl) \in L^{-\infty}(X^\wedge;\gL)$. Let 
$\tilde\go_2 \prec\tilde\go_1 \prec\tilde\go_0$ be cut-off functions, 
\[ b_1(\gl):= \tilde\go_1(r\ladm)\,r^{m}\op_M(h)(\gl)\tilde\go_0(r\ladm)
  +(1-\tilde\go_1(r\ladm))\,b_0(\gl)(1-\tilde\go_2(r\ladm)). \]
Then \ref{param1.1}. is satisfied because of \eqref{roughparam}, and because
$\op_M(h)(\gl)-b_0(\gl)\in L^{-\infty}(X^\wedge;\gL)$. 
Further, \ref{param1.2}. is a consequence of the holomorphy of $h$ 
and the fact that $\tilde q(0,x,r\gr,\xi,0)$ is the symbol of a
local parametrix of $\sumk{m} a_k(0)(-r\partial_r)^k$. Finally, since
$b_1(\gl)$ is a parameter-dependent parametrix of $a_{\wedge}(\gl)$ in the
interior, and because its operator
norm in $\mathcal{L}(\mathcal{K}^{s,\gg-m},\mathcal{K}^{s+m,\gg})$ can be
estimated by semi-norms of its local symbols, we get
$b_1\in S_{\cl}^{-m,m}(\gL;\mathcal{K}^{s,\gg-m},\mathcal{K}^{s+m,\gg})$. 
Recall that $a_{\wedge}$ is itself an operator-valued symbol of order
$(m,m)$, and the local symbols of $b_1(\gl)$ can be expressed in terms of 
those of $a_{\wedge}(\gl)$. Hence $a_{\wedge}b_1-1\in
S_{\cl}^{0,m}(\gL;\mathcal{K}^{s,\gg-m},\mathcal{K}^{s,\gg-m})$ for every
$s\in\R$. Together with \ref{param1.1}. this implies \ref{param1.3}.
\hfill $\qed$

\medskip
{\bf Proof of Lemma \ref{param2}.} 
A necessary condition for $a_{\wedge}(b_1+v_1)-1$ to be of Green type
is that its conormal symbol vanishes. So, we need 
$\sM^0(a_{\wedge}(b_1+v_1))(z)=1$ for $z\in\gG_{\frac{n+1}{2}-\gg}$. 
Because the relation \eqref{comp-co} is also valid for cone pseudodifferential
operators and because $\sM^{-m}(b_1+v_1)(z)=h(0,z,0)+f(z)$, we want 
$f$ to satisfy
\[ \sM^m(a_{\wedge})(z-m)\sM^{-m}(b_1+v_1)(z)=
    \sM^{m}(A)(z-m)\big(h(0,z,0)+f(z)\big)=1. \] 
Observe that the ellipticity of $A-\gl$ implies the invertibility of  
$\sM^{m}(A)(z)$ on $\gG_{\frac{n+1}{2}-\gg}$ so that for 
$z\in\gG_{\frac{n+1}{2}-\gg+m}$ we may set
\[ f(z):=\sM^{m}(A)^{-1}(z-m)\big(1-\sM^{m}(A)(z-m)h(0,z,0)\big). \]
This function is clearly meromorphic with poles contained in the boundary 
spectrum of $A$. Now, for $|\gl|\ge 1$ the operators of multiplication by
$\go_0(r\ladm)$ and $\go_1(r\ladm)r^m$ are twisted homogeneous of
degree $(0,m)$ and $(-m,m)$, respectively. Moreover, $v_1(\gl)$ is smoothing
for every $\gl$. Thus  
\[ v_1\in S_{\cl}^{-m,m}(\gL;\mathcal{K}^{s,\gg-m}(X^\wedge),
   \mathcal{K}^{s',\gg}(X^\wedge)^{\gs})\;\text{ for all } s,s',\gs\in\R. \]
Finally, the family $a_{\wedge}(\gl)(b_1(\gl)+v_1(\gl))-1$ is of Green type 
because it is a pointwise smoothing operator-valued symbol with 
$\sM^0(a_{\wedge}(b_1+v_1)-1)=0$.
\hfill $\qed$

\medskip
{\bf Proof of Lemma \ref{param3}.} 
Due to Lemma \ref{param2} there is 
$g_r\in\edgeg{0,m}{\gg-m,\gg-m}_{\eps_r}$ such that $a_{\wedge}(b_1+v_1) -
1=g_r$. In a similar way, we get $(b_1+v_1)a_{\wedge}-
1=g_l\;$ for some $g_l\in \edgeg{0,m}{\gg,\gg}_{\eps_l}$.
Because $\sK^m(a_{\wedge})=\sK^m(a)$ is invertible for $\gl\in\gL$ 
sufficiently large, we obtain there the elementary relation
\[ \sK^m(a_{\wedge})^{-1}=\sK^{-m}(b_1+v_1)-\sK^{-m}((b_1+v_1) g_r)
   +\sK^0(g_l)\sK^m(a_{\wedge})^{-1}\sK^0(g_r). \]
Let $\chi\in C^\infty(\gL)$ be such that $\chi\equiv 0$ where 
$\sK^m(a_{\wedge})$ is not invertible, and $\chi\equiv 1$ for $\gl$ 
sufficiently large. Then $\chi\sK^m(a_{\wedge})^{-1}\in
C^{\infty}(\gL,\mathcal{L}(\mathcal{K}^{s-m,\gg-m},\mathcal{K}^{s,\gg}))$ 
is twisted homogeneous of degree $(-m,m)$ and so in
$S_{\cl}^{-m,m}(\gL;\mathcal{K}^{s-m,\gg-m},\mathcal{K}^{s,\gg})$, cf.
Example~\ref{opvs6}. Now, for $g_0:=g_l\,\chi\sK^m(a_{\wedge})^{-1}g_r$,
Proposition~\ref{opvs3} yields
\begin{align*} 
 g_0   &\in S_{\cl}^{-m,m}\big(\gL;\mathcal{K}^{s,\gg-m}(X^\wedge),
        \mathcal{K}^{s',\gg+\eps_l}(X^\wedge)^{\gs}\big)\;\text{ and}\\
 g_0^* &\in S_{\cl}^{-m,m}\big(\gL;\mathcal{K}^{s,-\gg}(X^\wedge),
        \mathcal{K}^{s',-\gg+m+\eps_r}(X^\wedge)^{\gs}\big)
\end{align*}
for all $s,s',\gs\in\R$. Let $\eps=\min(\eps_l,\eps_r)$. 
The continuous embeddings 
\[ \mathcal{K}^{s',\gg+\eps_l}(X^\wedge)^{\gs}\emb
   \mathcal{K}^{s',\gg+\eps}(X^\wedge)^{\gs} \;\text{ and }\;
   \mathcal{K}^{s',-\gg+m+\eps_r}(X^\wedge)^{\gs}\emb
   \mathcal{K}^{s',-\gg+m+\eps}(X^\wedge)^{\gs} \]
together with Proposition~\ref{opvs3} imply
$g_0\in\edgeg{-m,m}{\gg-m,\gg}_\eps$. 
Thus the proof is done by setting $g_1:=-(b_1+v_1) g_r+g_0$. 
\hfill $\qed$

\medskip
{\bf Proof of Lemma \ref{param4}.} 
For any $N\in\N$ we have
\[ a_{\wedge}b_{\wedge}^0\sumj{N-1}(1-a_{\wedge}b_{\wedge}^0)^j
    =1-(1-a_{\wedge}b_{\wedge}^0)^N. \]
Lemma~\ref{param3} implies
$(1-a_{\wedge}b_{\wedge}^0)^N\in\edgeg{-N,m}{\gg-m,\gg-m}$ so that 
\begin{equation}\label{ed-p1}
 b_{\wedge}^{(N)}:=b_{\wedge}^0\sumj{N-1}(1-a_{\wedge}b_{\wedge}^0)^j
 = b_{\wedge}^0 + \sum_{j=1}^{N-1}b_{\wedge}^0(1-a_{\wedge}b_{\wedge}^0)^j
\end{equation}
is a rough parametrix of $a_{\wedge}$ with a better remainder than
$b_{\wedge}^0$. Moreover, we have $b_{\wedge}^0(1-a_{\wedge}b_{\wedge}^0)^j
\in\edgeg{-m-j,m}{\gg-m,\gg}$ for every $j\in\N$. Therefore, there is a 
Green symbol $g_2\in\edgeg{-m-1,m}{\gg-m,\gg}$ such that 
$g_2\sim\sum_{j=1}^\infty b_{\wedge}^0(1-a_{\wedge}b_{\wedge}^0)^j$.
Finally, 
\begin{equation}\label{ed-p2} 
 b_{\wedge}(\gl):= b_{\wedge}^0(\gl) +g_2(\gl)
 = b_1(\gl)+v_1(\gl)+g_1(\gl)+g_2(\gl) 
\end{equation}
is a parameter-dependent parametrix of $a_{\wedge}(\gl)$.
\hfill $\qed$

\medskip
{\bf Construction of a global parametrix}.
Let us write $a(\gl)=A-\gl$ as in \eqref{pdell4}, let 
$Q(\gl)\in L_{\cl}^{-m,m}(\widetilde{\mathbb{B}};\gL)$ be a parametrix of
$P-\gl$. With $b_{\wedge}$ from Lemma~\ref{param4} define
\begin{equation}\label{parametrix}
  b(\gl)= \go_1\,b_{\wedge}(\gl)\,\go_0+(1-\go_1)\,Q(\gl)\,(1-\go_2)
\end{equation}
for $\gl\in\gL$  and $\go_2\prec\go_1\prec\go_0$.
This parameter-dependent pseudodifferential operator of order $(-m,m)$ is a
compound of `local' parametrices. Next we want to show that $b(\gl)$ is in 
fact a global parametrix of $a(\gl)$. To this end let $\go$ and $\go_3$ be
cut-off functions such that $\go_0\prec\go$ and $\go_3\prec\go_2$. Then
\begin{align*}
b(\gl)a(\gl) &=\big[\go_1\,b_{\wedge}(\gl)\go_0
    +(1-\go_1)\,Q(\gl)(1-\go_2)\big]\,a(\gl)\\
  &=\go_1\,b_{\wedge}(\gl)\,\go_0\,a(\gl)\go 
    +(1-\go_1)\,Q(\gl)(1-\go_2)\,a(\gl)(1-\go_3)\\
  &\hspace*{1em}+\go_1\,b_{\wedge}(\gl)\big[\go_0\,a(\gl)(1-\go)\big]
    +(1-\go_1)\,Q(\gl)\big[(1-\go_2)\,a(\gl)\go_3\big]\\
  &=\go_1\,b_{\wedge}(\gl)\go_0\,a(\gl)\go
    +(1-\go_1)\,Q(\gl)(1-\go_2)\,a(\gl)(1-\go_3)\\ 
  &=\go_1\,b_{\wedge}(\gl)\go_0\,a_{\wedge}(\gl)\go
    +(1-\go_1)\,Q(\gl)(1-\go_2)(P-\gl)(1-\go_3)\\ 
  &=\go_1\,b_{\wedge}(\gl)a_{\wedge}(\gl)\go
    +(1-\go_1)\,Q(\gl)(P-\gl)(1-\go_3)\\
  &\hspace*{1em}-\big[\go_1\,b_{\wedge}(\gl)(1-\go_0)\big]a_{\wedge}(\gl)\go
    -\big[(1-\go_1)\,Q(\gl)\go_2\big](P-\gl)(1-\go_3). 
\end{align*}
Here we have used the local property of $a(\gl)$. Now, 
the terms inside the brackets in the last line are smoothing elements in 
their classes because $\supp(\go_1)\cap\supp(1-\go_0)$ and 
$\supp(1-\go_1)\cap\supp(\go_2)$ are both empty. 
Because of the inclusions 
\[ \edgeg{-\infty}{\gg,\gg}_\eps \subset
  \mathcal{S}(\gL,\mathcal{L}(\mathcal{K}^{s,\gg-m},\mathcal{K}^{s',\gg-m})) 
  \;\text{ and }\; L^{-\infty}(\widetilde{\mathbb{B}};\gL)
  \subset \mathcal{S}(\gL,\mathcal{L}(H^{s},H^{s'})) \]
for all $s$, $s'\in\R$, we can easily verify that 
$ba-1\in \mathcal{S}(\gL,\coneg{\gg,\gg}_\eps)$. Similarly, we obtain 
$ab-1\in \mathcal{S}(\gL,\coneg{\gg-m,\gg-m}_\eps)$ so Theorem~\ref{eparam} 
is proved. \hfill $\qed$

\begin{remark} \label{opv-symbol}
The given construction of $b_{\wedge}$ yields immediately
\[ b_{\wedge}\in S_{\cl}^{-m,m}(\gL;\mathcal{K}^{s-m,\gg-m}(X^\wedge),
  \mathcal{K}^{s,\gg}(X^\wedge))\;\text{ for all } s\in\R. \]
This implies that $b(\gl)$ from \eqref{parametrix} belongs to
$S_{\cl}^{-m+M,m}(\gL;\mathcal{H}^{s-m,\gg-m}(\mathbb{B}),
\mathcal{H}^{s,\gg}(\mathbb{B}))$ for some $M$ depending on 
the norm estimates of $Q(\lambda)$, and depending on the constants of growth
corresponding to $\norm{\gk(\gl)}_{\mathcal{L}(\mathcal{K}^{s,\gg})}$ and 
$\norm{\gk^{-1}(\gl)}_{\mathcal{L}(\mathcal{K}^{s-m,\gg-m})}$, 
cf. Lemma~\ref{opvs0}. The spaces on $\mathbb{B}$ are here equipped with 
the trivial group action $\kappa=\mathrm{id}$.
\end{remark}
\begin{remark} \label{opvs-new}
The definition of $b_1(\gl)$ in the proof of Lemma~\ref{param1},
by means of parameter-dependent cut-off functions, is just convenient for 
technical purposes. In that way it is easier to verify that $b_{\wedge}(\gl)$ 
is an operator-valued symbol. But in fact, there exists 
$g\in\edgeg{-m,m}{\gg-m,\gg}$ (cf. \cite[Theorem 3.18]{GSS2}) such that 
\begin{equation}\label{px.new}\tag{\ref{parametrix}$'$}
 b(\gl)= \go_1\big(r^{m}\op_M(h)(\gl)+ v(\gl)+g(\gl)\big)\go_0
         +(1-\go_1)\,Q(\gl)\,(1-\go_2)
\end{equation} 
is also a global parametrix of $A-\lambda$, where $h$ is the symbol found in
the proof of Lemma~\ref{param1}, $v(\gl)$ is like in Lemma~\ref{param2},
and $Q(\gl)$ is as in \eqref{parametrix}.
\end{remark} 

\subsection{Holomorphic weakly parametric symbols}\label{weak-exp}
Let $\gO\subset\R^n$ be an open set and $\gL$ be a sector in $\C$. 
Let $\gG_\gd=\{z\in\C\sod \Re z=\gd\}$.

\begin{definition}\label{p-hs0}
For $\mu\in\R$ define $S_{\mathcal{O}}^\mu(\gO\times\R^n\times\C)$ as the 
class of holomorphic functions $h\in\mathcal{O}(\C,S^\mu(\gO\times\R^n))$
such that $h|_{\gG_\gd}\in S^\mu(\gO\times\R^n;\gG_\gd)$ for each 
$\gd\in\R$, uniformly for $\gd$ in compact intervals. 

This is a Fr{\'e}chet space with the system of semi-norms given by
\begin{equation}\label{so-sn}
 \sup_{\gd\in I}\big\{ \sup_{\substack{(\xi,\gr)\in\R^{n+1}\\ x\in K}}
 \spk{\xi,\gr}^{-\mu+|\ga|+\ell}|\partial_x^\gb\partial_\xi^\ga
 \partial_\gr^\ell h(x,\xi,\gd+i\gr)|\big\}
\end{equation}
for $\ga,\gb\in\N_0^n$, $K\subset\!\subset\gO$, $\ell\in\N_0$ and
$I\subset\!\subset\R$.
\end{definition}

\begin{definition}\label{p-hs1}
For $d\in\N$ define $S_{\mathcal{O}}^{\mu,d}(\gO\times\R^n\times\C;\gL)$
as the class of parameter-dependent symbols consisting of all 
$h(x,\xi,z,\gl)$ such that
\begin{enumerate}
\item[(i)] $z\mapsto h(\cdot,\cdot,z,\cdot)\in
    \mathcal{O}(\C,S^{\mu,d}(\gO\times\R^n;\gL))$;
\item[(ii)] $h(\cdot,\cdot,\gd+i\gr,\gl)\in 
    S^{\mu,\bar d}(\gO\times\R^n;\gG_\gd\times\gL)$ for each $\gd\in\R$,
    uniformly for $\gd$ in compact intervals, and with $\bar d=(1,d)$
    \;(cf. Appendix~\ref{pdo-p}).
\end{enumerate}
\end{definition}
Furthermore, a symbol $h\in S_{\mathcal{O}}^{\mu,d}(\gO\times\R^n\times
\C;\gL)$ is said to depend holomorphically on $\gl$ if the map 
\begin{equation*}
 \gl\mapsto h(\cdot,\cdot,\cdot,\gl): 
 \gLint\to S^{\mu}(\gO\times\R^n;\gG_\gd)
\end{equation*}
is holomorphic for every $\gd\in\R$. Observe that (ii) implies 
$h(\cdot,\cdot,\cdot,\gl)\in S^{\mu}(\gO\times\R^n;\gG_\gd)$. 
Finally, if for every $\gd\in\R$ the parameter-dependent symbols in (ii) 
are asked to be classical, then we will write 
$h\in S_{\mathcal{O},\cl}^{\mu,d}(\gO\times\R^n\times\C;\gL)$. In the 
case when $h$ is classical there is, for every $\gd\in\R$, an expansion
\[ h\sim\sum h_{\mu-j}\quad\text{in }\;
   S^{\mu,\bar d}(\gO\times\R^n;\gG_\gd\times\gL) \]
such that for $|\xi|+|\gr|+|\gl|^{1/d}\geq 1$,
\[ h_{\mu-j}(x,\tau\xi,\gd +i\tau\gr,\tau^d\gl)=\tau^{\mu-j}
   h_{\mu-j}(x,\xi,\gd+i\gr,\gl) \;\text{ for every } \tau\geq 1. \]
When $\gd$ is fixed we will identify $\gG_\gd\cong\R$ and then replace
in the argument of the symbols $\gd+i\gr$ by $\gr$.

\begin{lemma}\label{p-hs2} 
If $h\in S_{\mathcal{O},\cl}^{\mu,d}(\gO\times\R^n\times\C;\gL)$ depends
holomorphically on $\gl$, then for any fixed $\gd\in\R$ every homogeneous
component $h_{\mu-j}\in S^{\mu-j,\bar d}(\gO\times\R^n;\gG_\gd\times\gL)$
is holomorphic on $\gLint\cap\{|\gl|>1\}$, as a function taking values
in $S^{\mu-j}(\gO\times\R^n;\gG_\gd)$.
\end{lemma}
\begin{proof}
For simplicity we omit the variable $x$ and assume $d=1$; the general case 
is completely analogous. We first show the holomorphy of the 
principal symbol $h_\mu$ making use of the relation 
\begin{equation}\label{p-limit}
 h_\mu(\xi,\gr,\gl)=\lim_{\tau\to\infty}\tau^{-\mu}h(\tau\xi,\tau\gr,\tau\gl).
\end{equation}
If the convergence in \eqref{p-limit} is compact, that is,
uniformly on any compact subset of $\gLint\cap\{|\gl|>1\}$, then $h_\mu$ is
holomorphic there since so is $h$; consequently, $h-h_\mu$ is also holomorphic
and we can proceed as above to assert the same for $h_{\mu-1}$. Induction then
yields the holomorphy of every $h_{\mu-j}$.

To prove the compact convergence of \eqref{p-limit} in the topology of 
$S^{\mu}(\R^n;\gG_\gd)$ let $\ga\in\N_0^n$ and $\ell\in\N_0$ be given. 
For $\tau\geq 1$ and $|\gl|\geq 1$ we have
\[ h_\mu(\xi,\gr,\gl)-\tau^{-\mu}h(\tau\xi,\tau\gr,\tau\gl)=
  \tau^{-\mu}\bigl(h_\mu(\tau\xi,\tau\gr,\tau\gl)-
  h(\tau\xi,\tau\gr,\tau\gl)\bigr). \]
Denoting $r_{\mu-1}:=h_\mu-h\in S_{\cl}^{\mu-1,\bar d}(\R^n;
\gG_\gd\times\gL)$ then
\begin{align*}
  \tau^{-\mu}\abs{\partial_\xi^\ga\partial_\gr^\ell
  r_{\mu-1}(\tau\xi,\tau\gr,\tau\gl)}
   &=\tau^{-\mu+|\ga|+\ell}\abs{(\partial_\xi^\ga\partial_\gr^\ell 
      r_{\mu-1})(\tau\xi,\tau\gr,\tau\gl)}\\
   &\leq C_{\ga \ell}\,\tau^{-\mu+|\ga|+\ell}
      \spk{\tau\xi,\tau\gr,\tau\gl}^{\mu-|\ga|-\ell}
      \spk{\tau\xi,\tau\gr,\tau\gl}^{-1}\\
   &\leq C_{\ga \ell}'\spk{\xi,\gr}^{\mu-|\ga|-\ell}
      \spk{\gl}^{|\mu-|\ga|-\ell|}\tau^{-1},
\end{align*}
because $\spk{\tau\eta,\tau\gl}\leq 2\tau\spk{\eta}\spk{\gl}$ and
$\spk{\tau\eta,\tau\gl}^{-1}\leq\tau^{-1}\spk{\eta}^{-1}$ for 
$|\gl|\geq 1$. So,
\[ \sup_{\xi,\gr}\spk{\xi,\gr}^{-\mu+|\ga|+\ell}
   \abs{\partial_\xi^\ga\partial_\gr^\ell \bigl(h_\mu(\xi,\gr,\gl)
   - \tau^{-\mu}h(\tau\xi,\tau\gr,\tau\gl)\bigr)}
   \leq C\spk{\gl}^{|\mu-|\ga|-\ell|}\tau^{-1}. \]
This implies that the convergence in \eqref{p-limit} is compact for $|\gl|>1$.
\end{proof}

{\bf Weakly parametric symbols.} 
For our purposes we need an anisotropic version of a symbol class
introduced by {\sc Grubb} and {\sc Seeley} in \cite{GrSe95}. 
Let $d$ denote the anisotropy. 

\begin{definition}\label{p-hs3}
For $d\in\N$ and $\nu\in\R$ define
$S_{\mathcal{W}}^{\nu,d}(\gO\times\R^n\times \R;\gL)$ as the space of
functions $h\in C^{\infty}(\gO\times\R^n\times\R \times\gL)$ that are
holomorphic in $\gl$ for $|\gl|>1$, and such that for every $k$,
\[ \partial_w^k h(\cdot,\cdot,\cdot,\tfrac{1}{w^d})\in
   S^{\nu+k}(\gO\times\R^n\times\R)\;\text{ for } \frac{1}{w^d}\in\gL,
   \;\text{ uniformly for } |w|\leq 1.\]
\end{definition}
Note that there are $d$ convex subsets of 
$\{w\in\C\sod \frac{1}{w^d}\in\gL\text{ and }|w|\leq 1\}$ such that 
$w\mapsto\frac{1}{w^d}$ maps each one onto $\gL\cap\{|\gl|\geq 1\}$.
We fix one of these subsets and denote it by $D(\gL)$.
The requirement for $h$ means that for every $k\in\N_0$,
\begin{equation}\label{pi-sn}
  \pi_{w,k}(h):=\sup_{\substack{(\xi,\gr)\in\R^{n+1}\\ x\in K}}
  \spk{\xi,\gr}^{-\nu-k+|\ga|+\ell}
  \abs{\partial_x^\gb\partial_\xi^\ga\partial_\gr^\ell\partial_w^k 
       h(x,\xi,\gr,\tfrac{1}{w^d})} 
\end{equation}
must be uniformly bounded in $D(\gL)$ for every $\ga,\gb\in\N_0^n$,
$K\subset\!\subset\gO$ and $\ell\in\N_0$.

\begin{theorem}\label{p-hs6}
Let $h\in S_{\mathcal{W}}^{\nu,d}(\gO\times\R^n\times\R;\gL)$, $d\in\N$ and
$\nu\in\R$. For any $N\in\N$ there are symbols 
$h_k\in S^{\nu+k}(\gO\times\R^n\times\R)$ such that
\begin{equation}\label{p-exp}
 \gl^{N/d}\Big\{h(x,\xi,\gr,\gl)-\sumk{N-1}\gl^{-k/d}h_k(x,\xi,\gr)\Big\}
 \in S_{\mathcal{W}}^{\nu+N,d}(\gO\times\R^n\times\R;\gL).
\end{equation}
\end{theorem}
\begin{proof}
The variable $x$ will be omitted again since it does not play any role along
the proof. Let us set
\[ f(\xi,\gr,w):=h(\xi,\gr,\tfrac{1}{w^d})\;\text{ for } w\in D(\gL). \]
For every $k\in\N_0$ the mapping $w\mapsto\partial_w^{k+1} f:D(\gL)\to 
S^{\nu+k+1}$ is bounded by definition, thus 
$w\mapsto\partial_w^{k} f:D(\gL)\to S^{\nu+k+1}$ is uniformly continuous
(fundamental theorem of calculus) implying that 
$h_k:=\frac{1}{k!}\lim_{w\to 0}\partial_w^k f(w)$ exists 
in $S^{\nu+k+1}(\R^n\times\R)$; note that $0\in\partial D(\gL)$. 
In particular, for every $(\xi,\gr)\in\R^n\times\R$ and $\ga\in\N_0^{n+1}$
we have the pointwise convergence $\partial_{\xi,\gr}^\ga\partial_w^k
f(\xi,\gr,w)\to \partial_{\xi,\gr}^\ga h_k(\xi,\gr)$ as $w\to 0$, and so
\[ \spk{\xi,\gr}^{-\nu-k+|\ga|}\partial_{\xi,\gr}^\ga\partial_w^k 
   f(\xi,\gr,w)\to \spk{\xi,\gr}^{-\nu-k+|\ga|}\partial_{\xi,\gr}^\ga
   h_k(\xi,\gr)\;\text{ as }\; w\to 0. \] 
Therefore, $h_k\in S^{\nu+k}(\R^n\times\R)$ since $w\mapsto\partial_w^{k}f$
is bounded in $S^{\nu+k}(\R^n\times\R)$.
Now, for every $N\in\N$ a Taylor expansion yields 
\[ w^{-N}\Bigl\{f(\xi,\gr,w)-\sumk{N-1}w^{k}h_k(\xi,\gr)\Bigr\}
  =\tfrac{1}{(N-1)!}\int_0^1(1-t)^{N-1}(\partial_w^N f)(\xi,\gr,tw)dt. \]
The left side, written in terms of $\gl=\tfrac{1}{w^d}$, is exactly the
expression in \eqref{p-exp}. Denoting the integral above by 
$r_N(\xi,\gr,\tfrac{1}{w^d})$ we have for every $k$,
\[ \partial_w^k r_N(\xi,\gr,{\tfrac{1}{w^d}})
  =\int_0^1(1-t)^{N-1}t^k(\partial_w^{N+k}f)(\xi,\gr,tw)dt. \]
We want to show that $\partial_w^k r_N(\cdot,\cdot,\tfrac{1}{w^d})
\in S^{\nu+N+k}(\R^n\times\R)$ uniformly for $w\in D(\gL)$, 
but this follows from the estimates
\begin{align*}
 \sup_{w\in D(\gL)} \pi_{w,k}(r_N)
   &\leq\int_0^1(1-t)^{N-1}t^k\sup_{w\in D(\gL)}\pi_{tw,N+k}(h)dt\\  
   &\leq\int_0^1(1-t)^{N-1}t^k\sup_{tw\in D(\gL)}\pi_{tw,N+k}(h)dt
\end{align*}  
taking into account that $h\in S_{\mathcal{W}}^{\nu,d}$ and so, in particular,
the supremum in the latter integral is uniformly bounded.
\end{proof}

The expansion \eqref{p-exp} is the main motivation for the consideration of 
the class $S_{\mathcal{W}}^{\nu,d}$. As it was done in \cite{GrSe95}, this 
kind of expansions can be used to obtain a complete asymptotic expansion of 
the resolvent of certain pseudodifferential operators on smooth compact
manifolds (cf. \cite{Ag87}). In particular, Dirac-type operators with nonlocal
boundary conditions in the spirit of {\sc Atiyah, Patodi} and {\sc Singer}
\cite{APS} are considered in \cite{GrSe95}. Further related results can be 
found in \cite{Gru99}, \cite{GrSc99} and \cite{GrSe96}.

\begin{lemma}\label{p-hs4}
Let $\nu\leq 0$ and $\gd\in\R$. Let $h\in S^{\nu,\bar d}
(\gO\times\R^n;\gG_\gd\times\gL)$ be homogeneous of degree $(\nu,\bar d)$ 
for $|\xi|+|\gr|+|\gl|^{1/d}\geq 1$, and such that it depends holomorphically
on $\gl$ for $|\gl|>1$. Then $h\in S_{\mathcal{W}}^{\nu,d}(\gO\times\R^n
\times\R;\gL)$, identifying $\gG_\gd\cong\R$.
\end{lemma}
\begin{proof} 
Without loss of generality we omit the variable $x$. It is true that
\[ h(\xi,\gr,\tfrac{1}{w^d})=|w|^{-\nu}h(|w|\xi,|w|\gr,|w|^d\tfrac{1}{w^d})
  \;\text{ for } w\in D(\gL).\]
Let $w=re^{i\gt}\in D(\gL)$, i.e., $r\leq 1$ and $e^{-id\gt}\in\gL$. Hence
\[ h(\xi,\gr,\tfrac{1}{w^d})=r^{-\nu}h(r\xi,r\gr,e^{-id\gt})=:f(\xi,\gr,r) \]
leading to the relation (because $h$ is holomorphic in $w$)
\[ \partial_w^k h(\xi,\gr,\tfrac{1}{w^d}) 
   =(\partial_r^k f)(\xi,\gr,r)e^{-ik\gt}\;\text{ for every } k\in\N. \]
To prove the assertion we need $\sup_{w\in D(\gL)}\pi_{w,k}(h)<\infty$
for every semi-norm as in \eqref{pi-sn}. Now, setting 
$\gl_0= e^{-id\gt}$ we have
\begin{align*}
(\partial_r^k f)(\xi,\gr,r) &=\partial_r^k\big(r^{-\nu}h(r\xi,r\gr,\gl_0)\big)
    =\sum_{\substack{j+j'=k\\ j'\leq-\nu}}
     C_{jj'}\,r^{-\nu-j'}\partial_r^j h(r\xi,r\gr,\gl_0)\\
  &=\sum_{\substack{j+j'=k\\ j'\leq-\nu}}
     C_{jj'}\,r^{-\nu-j'}\sum_{|\gb|+k'=j}C_{\gb k'}\,\xi^\gb\gr^{k'}
     (\partial_\xi^\gb\partial_\gr^{k'}h)(r\xi,r\gr,\gl_0).
\end{align*}
Because $|\xi^\gb\gr^{k'}|\leq\spk{\xi,\gr}^j$, and since 
$\partial_\xi^\gb\partial_\gr^{k'}h$ is a symbol of order $\nu-j$, then
for $\ga\in\N_0^n$ and $\ell\in\N$ the derivatives
$\partial_\xi^\ga\partial_\gr^\ell\partial_w^k h(\xi,\gr,\tfrac{1}{w^d})$ 
can be estimated by terms of the form
\begin{gather*}
 \sum_{\substack{j+j'=k\\ \nu+j'\leq 0}} C_{jj'}\,r^{-\nu-j'}\!\!
   \sum_{\ga_1+\ga_2=\ga}\;\sum_{\ell_1+\ell_2=\ell} C_{\ga \ell}
   \spk{\xi,\gr}^{j-|\ga_1|-\ell_1}
   \spk{r\xi,r\gr}^{\nu-j-|\ga_2|-\ell_2}\,r^{|\ga_2|+\ell_2}\\
 \qquad\leq\sum_{\substack{j+j'=k\\ \nu+j'\leq 0}} \tilde C_{jj'}\!\!
   \sum_{\ga_1+\ga_2=\ga}\;\sum_{\ell_1+\ell_2=\ell}
   \spk{\xi,\gr}^{j-|\ga_1|-\ell_1}\,r^{-\nu-j'+|\ga_2|+\ell_2}
   \spk{r\xi,r\gr}^{\nu+j'-|\ga_2|-\ell_2}.
\end{gather*}
Using now the inequality $r^{s}\spk{r\xi,r\gr}^{-s}\leq\spk{\xi,\gr}^{-s}$
for $r\leq 1$ and $s\geq 0$, we get
\[ \abs{\partial_\xi^\ga\partial_\gr^\ell\partial_w^k 
   h(\xi,\gr,\tfrac{1}{w^d})} \leq C\spk{\xi,\gr}^{\nu+k-|\ga|-\ell} \]
with a constant $C$ depending uniformly on $\gt=\arg{w}$ and being 
independent of the variable $r$. Hence $\,\sup_{w\in D(\gL)}\pi_{w,k}(h)$ 
is finite.
\end{proof}

\begin{lemma}\label{p-hs5}
Let $\nu\leq 0$ and $h\in S_{\mathcal{O},\cl}^{\nu,d}(\gO\times\R^n
\times\C;\gL)$. If $h$ depends holomorphically on $\gl$, its restriction 
from $\C$ to any $\gG_\gd$ induces a symbol in
$S_{\mathcal{W}}^{\nu,d}(\gO\times\R^n\times\R;\gL)$.
\end{lemma}
\begin{proof}
As it was done before we ignore for a moment the variable $x$. For any 
$\gd\in\R$ the symbol $h\in S_{\mathcal{O},\cl}^{\nu,d}(\R^n\times\C;\gL)$
admits an asymptotic expansion 
\[ r_N=h-\sumj{N-1}h_{\nu-j}\in S_{\cl}^{\nu-N,\bar d}(\R^n;\gG_\gd\times\gL)
   \;\text{ for any } N\in\N \]
with symbols $h_{\nu-j}\in S^{\nu-j,\bar d}$ that are homogeneous for 
$|\xi|+|\gr|+|\gl|^{1/d}\geq 1$, and are holomorphic in $\gl$ for 
$|\gl|>1$ due to Lemma~\ref{p-hs2}. In order to prove the claim we have to 
investigate $\sup_{w}\pi_{w,k}(h)$ for any semi-norm in 
$S_{\mathcal{W}}^{\nu,d}(\R^n\times\R;\gL)$.
To this end let $k\in\N_0$ be given. For any homogeneous component each
semi-norm is bounded due to Lemma~\ref{p-hs4}, so we only need to show
$\sup_{w}\pi_{w,k}(r_N)<\infty$ for some $N$. In fact, choosing $N\geq\nu+k$ 
we achieve
\begin{align*}
\Bigl|\partial_\xi^\ga\partial_\gr^\ell\partial_w^k
    r_N(\xi,\gr,\tfrac{1}{w^d})\Bigr|
 &\leq \sum_{k'=1}^k C_{k'}\,|w|^{-k-dk'}\,  
    |(\partial_\xi^\ga\partial_\gr^{\ell}
    \partial_\gl^{k'}r_N)(\xi,\gr,\tfrac{1}{w^d})|\\
 &\leq \sum C_{k'}\,|w|^{-k-dk'}
    \spk{\xi,\gr,\tfrac{1}{w}}^{\nu-N-|\ga|-\ell-dk'}\\
 &=\sum C_{k'}\,|w|^{-k-dk'}\spk{\xi,\gr,\tfrac{1}{w}}^{-k-dk'}
    \spk{\xi,\gr,\tfrac{1}{w}}^{\nu-N+k-|\ga|-\ell}\\
 &\leq C\spk{\xi,\gr}^{\nu+k-|\ga|-\ell}\text{ with $C>0$ independent of $w$}.
\end{align*}
In the last estimate one uses the relations
$|w|^{-s}\spk{\xi,\gr,\tfrac{1}{w}}^{-s}\leq C_s$ for $s\geq 0$, $|w|\leq 1$,
and $\spk{\xi,\gr,\tfrac{1}{w}}^{\nu-N+k}\leq\spk{\xi,\gr}^{\nu+k}$
for $N\geq\nu+k$. 
\end{proof}

\begin{theorem}\label{p-hs7}
For $\nu\leq 0$ and $\mu\in\N_0$ every symbol 
$h\in S_{\mathcal{O},\cl}^{\nu-\mu,d}(\gO\times\R^n\times\C;\gL)$ 
depending holomorphically on $\gl$ admits, for any $N\in\N$, the expansion 
\begin{equation}\label{p-expL}
 \gl^{(N+\mu)/d}\Big\{h(x,\xi,\gd+i\gr,\gl)-\sumk{N-1}\gl^{(-k-\mu)/d}
 h_k(x,\xi,\gd+i\gr)\Big\}\in S_{\mathcal{W}}^{\nu+N,d}
\end{equation}
with $h_k\in S_{\mathcal{O}}^{\nu+k}(\gO\times\R^n\times\C)$ given by
\begin{equation}\label{p-coef}
  h_k(x,\xi,z)=\frac{1}{k!}\lim_{w\to 0}
 \partial_w^k\Bigl(w^{-\mu} h(x,\xi,z,\tfrac{1}{w^d})\Bigr),\quad w\in D(\gL).
\end{equation}
\end{theorem}
\begin{proof}
Due to the inclusion $\gl^{\mu/d}S_{\mathcal{O},\cl}^{\nu-\mu,d}\subset
S_{\mathcal{O},\cl}^{\nu,d}$ it is enough to prove the statement for $\mu=0$.
For simplicity, let us drop the variable $x$. Lemma~\ref{p-hs5}
assures the existence of an expansion \eqref{p-exp} of $h|_{z\in\gG_\gd}$ 
for every $\gd\in\R$. As proven in Theorem~\ref{p-hs6} the 
coefficients $h_k$ are indeed given by \eqref{p-coef} on every line $\gG_\gd$.
It only remains to prove that $h_k\in S_{\mathcal{O}}^{\nu+k}
(\R^n\times\C)$. In other words, we have to verify:
\begin{enumerate}
\item[(i)] 
   $h_k\in S^{\nu+k}(\R^n\times\gG_\gd)$ uniformly for $\gd$ in 
   compact intervals, 
\item[(ii)] 
   $z\mapsto h_k(\cdot,z)\in\mathcal{O}(\C,S^{\nu+k}(\R^n))$.
\end{enumerate}
To this end let $I\subset\!\subset\R$ be a compact interval. Since every
semi-norm $\pi_{w,k}\bigl(h|_{\gG_\gd}\bigr)$ is uniformly bounded for 
$w\in D(\gL)$ and $\gd\in I$ then (i) follows by a pointwise consideration
as in the proof of Theorem~\ref{p-hs6}. Furthermore, the convergence in
\eqref{p-coef}, taking place in $S^{\nu+k+1}$, is uniform for $\gd\in I$ 
hence compact in $\C$. Thus $h_k$ is holomorphic with values in
$S^{\nu+k+1}(\R^n)$. By using now (i) and the fundamental theorem of 
calculus we obtain for neighboring $z$ and $\tilde z$ the estimates
\begin{align*} 
 \abs{\partial_{\xi}^\ga(h_k(\xi,z)-h_k(\xi,\tilde z))} 
   &\leq |z-\tilde z|\int_0^1\abs{\partial_{\xi}^\ga\partial_\gr 
         h_k(\xi,\tilde z+t(z-\tilde z))}dt \\
   &\leq |z-\tilde z|\;C\spk{\xi}^{\nu+k-|\ga|}
         \int_0^1\spk{\gr_t}^{|\nu+k-|\ga||}dt 
\end{align*}
with $\Im\tilde z\leq\gr_t\leq\Im z=\gr$. That is, the mapping
$z\mapsto h_k(\cdot,z):\C\to S^{\nu+k}(\R^n)$ is continuous which
implies, in particular, that the Cauchy integral 
$\frac{1}{2\pi i}\int_{\partial B}\frac{h_k(\xi,\zeta)}{\zeta-z}d\zeta$ 
exists in $S^{\nu+k}(\R^n)$ for every disk $B$. 
Hence (ii) holds and the proof is done.
\end{proof}

\begin{remark}
This theorem will be applied to the homogeneous components of the 
pseudodifferential Mellin symbols appearing in the parametrix construction
of $A-\gl$, cf. Lemma~\ref{param1}. These homogeneous symbols clearly
depend holomorphically on $\gl$.
\end{remark}

\section{Heat trace asymptotics}
Along this section $A$ will always denote an elliptic cone differential
operator of order $m>0$ on the manifold $\mathbb{B}$, 
cf. Section~\ref{cone-do}.

\subsection{The operator $\mathbf{e^{-tA}}$ in the cone algebra}
\label{heat-op}

Fix $\delta>0$ and $0<\gp<\frac\pi{2}$. We consider the contour 
$\gU=\gU_1\cup\gU_2\cup\gU_3$ in $\C$, where
\begin{alignat*}{2}
\gl &=re^{i\gp} && (+\infty>r\geq \delta)\;\text{ on }\gU_1,\\
\gl &=\delta e^{i\gt} && (\gp\leq\gt\leq 2\pi-\gp)\;\text{ on }\gU_2,\\ 
\gl &=re^{-i\gp} &\quad & (\delta\leq r<+\infty)\;\text{ on }\gU_3.
\end{alignat*}

Let $\gL$ be the sector $\{\gl\in\C\sod \gp\le \arg\lambda\le 2\pi-\gp\}$
and let $\gL_\delta=\gL\cap\{|\lambda|\ge\delta\}$. We assume that $A-\gl$ is
invertible on an open neighborhood $\mathcal U(\gL_\delta)$ of $\gL_\delta$.
Using the identity $(A-\lambda)^{-1}\big(1-(\lambda-\lambda_0)
(A-\lambda_0)^{-1}\big)=(A-\lambda_0)^{-1}$ and the embedding properties of 
our weighted Sobolev spaces, it can be easily proven (via Neumann series)
that the map $\lambda\mapsto(A-\gl)^{-1}: \mathcal U(\gL_\delta)\to
\mathcal{L}(\mathcal{H}^{s-m,\gg-m}(\mathbb{B}),
\mathcal{H}^{s,\gg}(\mathbb{B}))$ is holomorphic. 

For $t>0$ we now define, as usual, 
\begin{equation}\label{etA}
  \etA =\frac{i}{2\pi}\intU e^{-t\gl}(A-\gl)^{-1}d\gl. 
\end{equation}
This integral converges absolutely in
$\mathcal{L}(\mathcal{H}^{s-m,\gg-m}(\mathbb{B}),
\mathcal{H}^{s,\gg}(\mathbb{B}))$ for any real $s$, and it does not depend on
$\gU$ provided that $\gU$ is a path in $\gL_\delta$ such that 
$\Re\gl\to\infty$ as $|\gl|\to\infty$. In particular, for $0<t<1$ the path
$t^{-1}\gU$ stays inside $\gL_\delta$ and we can write  
\begin{equation}\label{ettA} \tag{\ref{etA}$'$}
  \etA =\frac{i}{2\pi}\inttU e^{-t\gl}(A-\gl)^{-1}d\gl. 
\end{equation}
This representation will be more convenient for some computations later on.

\begin{theorem}\label{heat1} 
Let $A-\gl$ be parameter-elliptic with respect to $\gg$. Then, for every
$t>0$ the operator $\etA$ belongs to $\coneg{\gg-m,\gg}$.
\end{theorem}
\begin{proof}
Let $t>0$ be fixed. We need to show certain mapping properties according to
Definition~\ref{green1}. First of all, observe that (integration by parts)
\begin{equation}\label{etAl} \tag{\ref{etA}$''$}
  \etA =\tfrac{i(\ell-1)!}{2\pi}\, t^{-\ell+1}
  \intU e^{-t\gl}(A-\gl)^{-\ell}d\gl 
\end{equation}
for every $\ell\in\N$. Thus $\etA$ and its formal
adjoint $(\etA)^*$ are for every $s,s'\in\R$ bounded operators in
$\mathcal{L}(\mathcal{H}^{s,\gg-m},\mathcal{H}^{s',\gg})$ and
$\mathcal{L}(\mathcal{H}^{s,-\gg},\mathcal{H}^{s',-\gg+m})$, respectively.
Hence they are smoothing operators in the interior, but in order to be in
$\coneg{\gg-m,\gg}$ these operators have to improve
the weights by some $\eps>0$. We prove this by making use of the parametrix
$b(\lambda)$ of $A-\lambda$ given in \eqref{px.new}. Then $b(\lambda)^\ell$
is a parametrix of $(A-\gl)^{\ell}$ and $(A-\gl)^{-\ell} - b(\lambda)^\ell$
belongs to $\mathcal{S}(\gL_\delta,\coneg{\gg-m,\gg})$. Moreover, for 
$\ell>1$ we can write
\begin{equation}\label{l-parametrix}
 b(\lambda)^\ell= \go_1\left(r^{m\ell}\op_M(h_{\ell})(\gl) +
 g_\ell(\lambda)\right)\go_0 +(1-\go_1)\,Q_\ell(\gl)\,(1-\go_2) +
 G_\ell(\lambda),
\end{equation}
where 
$h_\ell\in C^\infty(\R_+,M_\mathcal{O}^{-m\ell,m}(X;\gL)$,
$g_\ell\in\edgeg{-m\ell,m}{\gg-m,\gg}$, 
$Q_\ell\in L_{\cl}^{-m\ell,m}(\widetilde{\mathbb{B}};\gL)$, and
$G_\ell\in \mathcal{S}(\gL,\coneg{\gg-m,\gg}$. 
Notice that all the contributions involving the smoothing Mellin element
$v(\gl)$ from \eqref{px.new} are now contained in $g_\ell(\lambda)$
or in $G_\ell(\lambda)$ since they are either supported in the interior of
$\mathbb{B}$ or multiplied by a Green element or multiplied by a factor
$r^{m\ell}$ which improves the weight whenever $\ell>1$. Clearly, the
nonsmoothing parts of $b(\lambda)^\ell$ improve the weight $\gg$ at least by
$m(\ell-1)$ while the families $g_\ell(\gl)$ and $G_\ell(\lambda)$ has the 
same gain $\eps>0$ as the smoothing part of $b(\gl)$.
\end{proof}

\subsection{Approximation of the resolvent}\label{hop-as}
Following {\sc Seeley}'s ideas \cite{See67} we want to approximate the
resolvent $(A-\gl)^{-1}$ by means of a suitable parameter-dependent 
parametrix of $A-\gl$. More generally, we will approximate $(A-\gl)^{-\ell}$
for any $\ell\in\N$. Our aim is to extract a finite part of the parametrix
having homogeneous components and such that the remainder decreases fast
enough in $\gl$. To this end, we make use of the operator-valued symbolic
calculus. For simplicity, we will mostly work with the case $\ell=1$ omitting
it from the notation; corresponding comments about the general case will be
given when necessary. 

As given in \eqref{px.new} consider the parametrix 
\[ b(\gl)= \go_1 \big(r^{m}\op_M(h)(\gl)+ v(\gl)+g(\gl)\big)\go_0 
  + (1-\go_1) Q(\gl) (1-\go_2) \]
of $A-\gl$, where $h(r,z,\gl)=\tilde h(r,z,r^m\gl)$ with $\tilde h\in 
C^{\infty}(\rpbar,M_\mathcal{O}^{-m,m}(X;\gL))$, $v(\gl)$ is a family of
smoothing operators as in Lemma~\ref{param2}, and
$g\in\edgeg{-m,m}{\gg-m,\gg}$. 
Recall that $Q(\gl)$ is defined as an asymptotic summation
$Q(\gl)\sim\sum_j Q_j(\gl)$, where each $Q_j(\gl)$ is a parameter-dependent 
pseudo\-differential operator of order $-m-j$ with anisotropic homogeneous
local symbols. Moreover, since $v+g$ is a classical operator-valued symbol
of order $(-m,m)$, it admits an expansion of the form
\begin{equation} \label{r-e5}
  v(\gl)+g(\gl)=\sum_{j=0}^{N-1} g_{j}(\gl)+g_{[N]}(\gl)
\end{equation}
such that every $g_{j}$ is twisted homogeneous of degree $(-m-j,m)$, and
\[ g_{[N]}\in S^{-m-N,m}\big(\gL;
 \mathcal{K}^{s,\gg-m}(X^\wedge),\mathcal{K}^{s',\gg}(X^\wedge)\big)
 \text{ for every } s, s'\in\R. \] 
Thus the parametrix $b(\gl)$ can be expanded as 
\begin{equation} \label{r-e4.0}
\begin{split}
 b(\gl) &= b^{(N)}(\gl) + r_{N}(\gl) \\
 &:=\go_1\, b_{\wedge}^{(N)}(\gl)\, \go_0 +
  (1-\go_1) Q^{(N)}(\gl)(1-\go_2) + r_N(\gl),
\end{split}
\end{equation}
where $r_N(\gl)$ is a remainder of order $-m-N$, and 
\begin{equation*} 
 b_{\wedge}^{(N)}(\gl)= r^m\op_M(h)(\gl) + \sum_{j=0}^{N-1} g_j(\gl) 
 \;\text{ and }\; Q^{(N)}(\gl) = \sum_{j=0}^{N-1} Q_j(\gl).
\end{equation*}
Now, by means of the Taylor expansion 
\begin{align*}
 \tilde h(r,z,\gl) 
 &= \sumj{N-1} \frac{1}{j!}(\partial_r^j\tilde h)(0,z,\gl)r^j
    + r^N\tilde h_{{[N]}}(r,z,\gl) \\
 \text{with}\quad \tilde h_{{[N]}}(r,z,\gl) 
 &=\frac{1}{(N-1)!} \int_0^1(1-\gt)^{N-1}
    (\partial_r^N\tilde h)(\gt r,z,\gl)d\gt.
\end{align*}
Set $h_j(r,z,\gl):=(\partial_r^j\tilde h)(0,z,r^m\gl)$ for $j<N$, 
$h_{{[N]}}(r,z,\gl):=\tilde h_{{[N]}}(r,z,r^m\gl)$, and define 
\begin{equation} \label{r-e8}  
 b_{{\wedge},j}(\gl):= r^{m+j}\op_M(h_j)(\gl) + g_j(\gl).
 \quad \text{ (cf. \eqref{r-e5})}
\end{equation}
Then every $b_{{\wedge},j}$ is twisted homogeneous of degree $(-m-j,m)$, 
and we have
\begin{equation} \label{r-e4.2} 
 b_{\wedge}^{(N)}(\gl)= \sum_{j=0}^{N-1} b_{{\wedge},j}(\gl) +
 r^{m+N}\op_M(h_{{[N]}})(\gl). 
\end{equation}

Similarly, for $\ell\in\N$, $b(\gl)^\ell$ admits an expansion of the form
\begin{equation}\label{r-e6}
 b(\gl)^\ell = b_{\ell}^{(N)}(\gl) + r_{\ell,N}(\gl) 
\end{equation}
which can be used in order to approximate $(A-\gl)^{-\ell}$. 

\begin{theorem}\label{r-as2}
Given $k\in\N_0$ there exist $N(k)\in\N$ and $C_{k}>0$ such that 
\begin{equation*}
  \norm{(A-\gl)^{-\ell}- b_{\ell}^{(N)}(\gl)}_{\mathcal{L}
  (\mathcal{H}^{-k,\gg-m}(\mathbb{B}),\mathcal{H}^{k,\gg}(\mathbb{B}))}
  \leq C_{k}(1+|\gl|)^{-(k+\ell)}
\end{equation*}
for every $N\ge N(k)$ and $\gl\in\gL_\delta$.
\end{theorem}
\begin{proof} 
First of all, we have
\begin{align*}
 (A-\gl)^{-\ell}-b_{\ell}^{(N)}(\gl)
 &=\big\{(A-\gl)^{-\ell}-b(\gl)^\ell\big\}
   +\big\{b(\gl)^\ell-b_{\ell}^{(N)}(\gl)\big\} \\
 &=\big\{(A-\gl)^{-\ell}-b(\gl)^\ell\big\}
   + r_{\ell,N}(\gl). 
\end{align*}
Since $(A-\gl)^{-\ell}-b(\gl)^{\ell}$ belongs to
$\mathcal{S}(\gL_\delta,\coneg{\gg-m,\gg})$, it is for every integer
$k$ an element of $\mathcal{L}(\mathcal{H}^{-k,\gg-m}(\mathbb{B}),
\mathcal{H}^{k,\gg}(\mathbb{B}))$, and its norm is $O(|\gl|^p)$ 
for every $p\in\mathbb{Z}$. Thus it remains to verify the norm estimate 
for $r_{\ell,N}(\gl)$. For simplicity of notation we will check it 
explicitly only for $\ell=1$. In this case, $r_{N}(\gl)$ is just
\[ b-b^{(N)}=\go_1\big(b_{\wedge}-b_{\wedge}^{(N)}\big)\go_0+
   (1-\go_1)\big(Q-Q^{(N)}\big)(1-\go_2), \]
where $b_{\wedge}(\gl)=r^{m}\op_M(h)(\gl)+ v(\gl)+g(\gl)$. From the 
standard parameter-dependent calculus (cf. \cite[Section 9]{Shubin}) 
we get the norm estimate
\begin{equation}\label{r-e2}
 \norm{(1-\go_1)\big(Q(\gl)-Q^{(N)}(\gl)\big)(1-\go_2)}_{\mathcal{L}
 (H^{-k}(\widetilde{\mathbb{B}}), H^{k}(\widetilde{\mathbb{B}}))} 
 \leq C_k'\,[\gl]_m^{-m-N+2k}
\end{equation}
for $N\geq-m+2k$. On the other hand, 
\[ \go_1\big(b_{\wedge}-b_{\wedge}^{(N)}\big)\go_0\in
   S_{\cl}^{-m-N,m}\big(\gL;\mathcal{K}^{-k,\gg-m}(X^\wedge),
   \mathcal{K}^{k,\gg}(X^\wedge)\big) \]
so that 
\begin{equation}\label{r-e3}
\norm{\go_1\big(b_{\wedge}(\gl)-b_{\wedge}^{(N)}(\gl)\big)\go_0}_{
\mathcal{L}(\mathcal{K}^{-k,\gg-m}(X^\wedge),\mathcal{K}^{k,\gg}(X^\wedge))} 
\leq C_k''\,[\gl]_m^{-m-N+M(k)}
\end{equation}
with $M(k)=M_1+M_2$, where $M_1$ and $M_2$ are the constants of growth 
corresponding to $\norm{\gk(\gl)}_{\mathcal{L}(\mathcal{K}^{k,\gg})}$ and 
$\norm{\gk^{-1}(\gl)}_{\mathcal{L}(\mathcal{K}^{-k,\gg-m})}$, 
respectively (cf. Lemma~\ref{opvs0}). 
By means of the estimates \eqref{r-e2} and \eqref{r-e3}
together with Lemma~\ref{op-norm} we get 
\[ \norm{b(\gl)-b^{(N)}(\gl)}_{\mathcal{L}
   (\mathcal{H}^{-k,\gg-m},\mathcal{H}^{k,\gg})}
   \leq \tilde C_k\,[\gl]_m^{-(k+1)m}\leq C_k(1+|\gl|)^{-(k+1)} \]
for $-m-N+\max(2k,M(k))\leq -(k+1)m$. Set $N(k)=km+\max(2k,M(k))$. 

For an arbitrary $\ell\in\N$ the proof is quite the same once 
$b(\gl)^\ell$ is split as $b(\gl)$ in \eqref{r-e4.0}. This can be done by
expanding $g_\ell(\gl)$ and $Q_\ell(\gl)$ in \eqref{l-parametrix} inside 
their corresponding pseudodifferential calculi. Notice that $G_\ell(\gl)$
in \eqref{l-parametrix} is rapidly decreasing.
\end{proof}

The finite parametrix $b_{\ell}^{(N)}(\gl)$ from \eqref{r-e6} can also 
be used to approximate the operator $\etA$ for small $t$. 
For $\ell, N\in\N$ and $0<t<1$ we define
\begin{equation} \label{h-e1}
 E_\ell^{(N)}(t):=
 \tfrac{i(\ell-1)!}{2\pi}\, t^{-\ell+1}\inttU e^{-t\gl}
 \,b_{\ell}^{(N)}(\gl)\,d\gl, 
\end{equation}
cf. \eqref{ettA} and \eqref{etAl}.
As a consequence of Theorem~\ref{r-as2} we get the following

\begin{corollary}\label{h-as3}
Given $k\in\N$ there exist $N(k)\in\N$ and $C_{k,\ell}>0$ such that
\[ \Big\|\etA-E_\ell^{(N)}(t)\Big\|_{\mathcal{L}
   (\mathcal{H}^{-k,\gg-m}(\mathbb{B}),
    \mathcal{H}^{k,\gg}(\mathbb{B}))}\leq C_{k,\ell}\,t^{k} \]
for every $N\geq N(k)$ and $0<t<1$.
\end{corollary}
\begin{proof} 
We combine \eqref{ettA} and \eqref{etAl} to write
\begin{equation*}
 \etA =\tfrac{i(\ell-1)!}{2\pi}\, t^{-\ell+1}
 \inttU e^{-t\gl}(A-\gl)^{-\ell}d\gl. 
\end{equation*}
Let $k\in\N$. In $\mathcal{L}(\mathcal{H}^{-k,\gg-m}(\mathbb{B}) 
\mathcal{H}^{k,\gg}(\mathbb{B}))$ we get
\begin{align*}
 \Big\|\etA-E_\ell^{(N)}(t)\Big\|
 &\le \tfrac{(\ell-1)!}{2\pi}\, t^{-\ell}\intU e^{-\Re\gl}
  \norm{(A-\tfrac\gl{t})^{-\ell}- b_{\ell}^{(N)}(\tfrac\gl{t})}|d\gl|\\
 &\le \tfrac{(\ell-1)!}{2\pi} C_{k}\, t^{-\ell} 
  \intU e^{-\Re\gl}(1+t^{-1}|\gl|)^{-(k+\ell)}|d\gl|\\
 &\le \tfrac{(\ell-1)!}{2\pi} C_{k}\, t^{k}
  \intU e^{-\Re\gl}|\gl|^{-(k+\ell)}|d\gl|\le C_{k,\ell}\,t^k
\end{align*}
after making the change of variables $\gl\to t^{-1}\gl$, and applying 
Theorem~\ref{r-as2}.
\end{proof}

\subsection{Asymptotic expansion of the heat trace}\label{trace-as}
As proved in Theorem~\ref{heat1} the operator $\etA$ belongs to
$\coneg{\gg-m,\gg}$ for $t>0$, so it is an operator of trace class
with kernel in $\mathcal{H}^{\infty,(-\gg+m+\eps)^{-}}(\mathbb{B})
\pitensor \mathcal{H}^{\infty,(\gg+\eps)^{-}}(\mathbb{B})$ for some
$\eps>0$, cf. Appendix~\ref{s-cone}. 
The kernel $K(t,y,y')$ of $\etA$ is commonly called the {\em heat kernel}. 
We further call the trace of $\etA$ the {\em heat trace} for $A$. 
In this section we will obtain an asymptotic expansion, as $t\to 0^+$, of 
the heat trace for a cone differential operator $A$ such that $A-\gl$ is
parameter-elliptic on the sector $\gL$ defined in Section~\ref{heat-op}. 
The complete expansion will be given in Theorem~\ref{htrace-exp}.

In order to expand the heat trace we will use the approximation
$E_\ell^{(N)}(t)$ from \eqref{h-e1}. Imitating the steps around
\eqref{r-e5}--\eqref{r-e4.2} for $b(\gl)^\ell$ we get in \eqref{r-e6} 
\begin{equation} \label{decomp1}
\begin{split}
 b_\ell^{(N)}(\gl) = 
 & \sum_{j=0}^{N-1} \left[\go_1\,b_{\wedge,\ell,j}(\gl)\,\go_0 + 
  (1-\go_1) Q_{\ell,j}(\gl)(1-\go_2)\right] \\
 & \qquad + \go_1\,r^{m\ell+N}\op_M(h_{\ell,[N]})(\gl)\,\go_0, 
\end{split}
\end{equation}
where every $Q_{\ell,j}$ is a parameter-dependent operator
of order $-m\ell-j$ with anisotropic homogeneous local symbols, 
$b_{{\wedge},\ell,j}=r^{m\ell+j}\op_M(h_{\ell,j})+g_{\ell,j}$ is 
twisted homogeneous of degree $(-m\ell-j,m)$, and $h_{\ell,[N]}$ is 
essentially the remainder of the Taylor expansion of the
parameter-dependent Mellin symbol $h_\ell$ from \eqref{l-parametrix}.

Using now \eqref{decomp1}, the integral \eqref{h-e1} becomes
\begin{equation*}
 E_\ell^{(N)}(t)=\sum_{j=0}^{N-1} E_{\ell,j}(t) + R_N(t) 
\end{equation*}
with
\begin{gather} \notag
 E_{\ell,j}(t) 
 = \tfrac{i(\ell-1)!}{2\pi}\, t^{-\ell+1} \inttU e^{-t\gl}
  \left[\go_1\,b_{\wedge,\ell,j}(\gl)\,\go_0 + 
  (1-\go_1) Q_{\ell,j}(\gl)(1-\go_2)\right] d\gl, \\ \label{R_N}
 R_N(t) = \tfrac{i(\ell-1)!}{2\pi}\, t^{-\ell+1}\inttU e^{-t\gl}\,
 \go_1\,r^{m\ell+N}\op_M(h_{\ell,[N]})(\gl)\,\go_0\, d\gl. 
\end{gather}
Every component supported in the interior of $\mathbb{B}$
\[ \tfrac{i(\ell-1)!}{2\pi}\, t^{-\ell+1} 
  \inttU e^{-t\gl}(1-\go_1) Q_{\ell,j}(\gl)(1-\go_2)\, d\gl \] 
can be treated on the closed manifold $\widetilde{\mathbb{B}}$ as in
\cite{Gilkey}. The components near the boundary
\begin{equation*} 
 \tfrac{i(\ell-1)!}{2\pi}\, t^{-\ell+1}
 \inttU e^{-t\gl}\,\go_1\,b_{\wedge,\ell,j}(\gl)\,\go_0\,d\gl  
\end{equation*}
and $R_N(t)$ require, however, more sophisticated calculations. 
More precisely, in order to expand in $t$ these boundary components, we will
use the twisted homogeneity of every $b_{\wedge,\ell,j}$, and the results 
from Section~\ref{weak-exp}. Notice that every $b_{{\wedge},\ell,j}$ 
is a sum of a parameter-dependent Mellin operator and an operator-valued 
Green symbol. While Mellin operators can be described through its local 
symbols, the operator-valued Green symbols are of global nature. 
For this reason, we will discuss these cases independently. 

If $a(\gl)$ is an operator family on $X^\wedge=\R_+\times X$ with Schwartz
kernels $k_a(\gl,r,x,r',x')$, then for $0<t<1$ and $x\in X$ we define formally 
\[ s[a](t,x):= t^{-\ell+1}\!\int_0^\infty\!\!    
   \inttU e^{-t\gl}\,\go_1(r)k_a(\gl,r,x,r,x)\dbar\gl dr \] 
with $\dbar\gl:=\frac{i}{2\pi}\, d\gl$. 
Without loss of generality we may assume $\go_1(r)=1$ for $r\le 1$, and
$\gd=1$ in the definition of $\gU$, cf. Section~\ref{heat-op}.

\bigskip
{\bf Expansion of the Green elements}
\begin{lemma}\label{tr-as3}
Let $g\in\edgeg{-m\ell-j,m}{\gg-m,\gg}$ be a twisted homogeneous Green
symbol of degree $(-m\ell-j,m)$, and let $k_g$ be its integral
kernel. Then for every $x\in X$, $0<t<1$ and $N\in\N$ with $N>j$, we have
\[ s[g](t,x)=t^{j/m}\left\{\int_0^\infty\!\!
   \intU e^{-\gl} k_g(\gl,r,x,r,x)\dbar\gl dr\right\}+O(t^{N/m}). \]   
\end{lemma}
\begin{proof}
Let $0<t<1$. Split the integral 
$\int_0^\infty=\int_0^{t^{1/m}}+\int_{t^{1/m}}^\infty$ in $s[g]$,
and denote the components by $s'[g]$ and $s''[g]$, respectively.
Because $\go_1(r)=1$ on $[0,t^{1/m}]$ we have
\[ s'[g](t,x)=t^{-\ell+1}\!\int_0^{t^{1/m}}\!\!
   \inttU\etgl k_g(\gl,r,x,r,x)\dbar\gl dr. \]
With the change $\gl\to t^{-1}\gl$ and Lemma~\ref{tr-as2} we get
\begin{align*}
\inttU\etgl k_g(\gl,r,x,r,x)d\gl
  &=t^{-1}\intU e^{-\gl} k_g(t^{-1}\gl,r,x,r,x)d\gl \\
  &=t^{-1+\ell+(j-1)/m}\!\intU e^{-\gl} k_g(\gl,t^{-1/m}r,x,t^{-1/m}r,x)d\gl.
\end{align*}
Therefore, the change of variables $r\to t^{1/m}r$ yields
\begin{align*}
s'[g](t,x) &=t^{(j-1)/m}\!\int_0^{t^{1/m}}\!\!
  \intU e^{-\gl} k_g(\gl,t^{-1/m}r,x,t^{-1/m}r,x)\dbar\gl dr\\
  &=t^{j/m}\left\{\int_0^{1}\!\!
  \intU e^{-\gl} k_g(\gl,r,x,r,x)\dbar\gl dr\right\}.
\end{align*}
With the same calculations we also get
\[ s''[g](t,x)=t^{j/m}\int_1^{\infty}\!\!
   \intU e^{-\gl}\go_1(t^{1/m}r) k_g(\gl,r,x,r,x)\dbar\gl dr. \]
A Taylor expansion at $r=0$ yields $\go_1(r)=1+r^M\go_{1,M}(r)$ with
\[ \go_{1,M}(r)=\frac{1}{(M-1)!}\int_0^1(1-\gt)^{M-1}\go_1^{(M)}(\gt r)\,d\gt
   \;\text{ and }\; \go_1^{(M)}=\partial_r^M\go_1.\]
In particular, for $M=N-j$ we get
\begin{equation}\label{tr-as3.1}
\begin{split}
s''[g](t,x) 
  &= t^{j/m}\left\{\int_1^{\infty}\!\!
  \intU e^{-\gl} k_g(\gl,r,x,r,x)\dbar\gl dr\right\}\\
  &\quad + t^{N/m}\int_1^{\infty}\!\!
  \intU e^{-\gl}\go_{1,N-j}(t^{1/m}r)\,r^{N-j} k_g(\gl,r,x,r,x)\dbar\gl dr.
\end{split}
\end{equation}
Since $\go_{1,M}(t^{1/m}r)$ is for any $M$ uniformly bounded, and because
$k_g$ decays rapidly as $r\to\infty$, the last integral is uniformly 
bounded as desired. 
\end{proof}

{\bf Expansion of parameter-dependent Mellin operators}
\medskip

First of all, we fix $\ell$ large enough such that every $E_{\ell,j}(t)$
and $R_N(t)$ from \eqref{R_N} have continuous Schwartz kernels.
For simplicity we will often drop $\ell$ from the notation.

Let $\tilde h\in M_{\mathcal{O}}^{-m\ell,m}(X;\gL)$ and set 
$h(r,z,\gl)=\tilde h(z,r^m\gl)$. Fix $j\in\N_0$ and define
\[ a(\gl):=r^{m\ell+j}\op_M(h)(\gl) \]
which is twisted homogeneous of degree $(-m\ell-j,m)$. Observe that the
Mellin component of $b_{\wedge,\ell,j}$ in \eqref{decomp1} is of this form.
Let $k_a(\gl,r,x,r',x')$ denote the continuous Schwartz kernel of $a(\gl)$,
taken with respect to the measure $dr'dx'$. As in Lemma~\ref{tr-as3} we
coveniently split the function $s[a](t,x)$ in two components
$s'[a](t,x)$ and $s''[a](t,x)$, and treat them separately.

\begin{lemma}\label{tr-as4}
For $0<t<1$ and $x\in X$ we have
\[ s'[a](t,x)=t^{j/m}\left\{\int_0^{1}\!\!\intU e^{-\gl}\,
   k_a(\gl,r,x,r,x)\dbar\gl dr\right\}. \]
\end{lemma}
\begin{proof}
Noting that $k_{a}(\gl,r,x,r,x)$ satisfies the relation \eqref{tr-a.e2} 
from Lemma~\ref{tr-as2} with $d=m$ and $\mu=-m\ell-j\;$ ($a(\gl)$ is twisted
homogeneous), the assertion follows making the same calculations as in the
proof of Lemma~\ref{tr-as3}.
\end{proof}

The second component of $s[a]$ is more delicate and cannot be treated as in
Lemma~\ref{tr-as3} since the kernels $k_a(\gl,r,x,r,x)$ increase as
$r\to\infty$, so the integrals without a cut-off function do not exist.
However, we can get an expansion in $t$ of $s''[a](t,x)$ by means of
the results from Section~\ref{weak-exp} (cf. Theorem~\ref{p-hs7}).
Recall that
\begin{equation*}
 s''[a](t,x)= t^{-\ell+1}\!\int_{t^{1/m}}^\infty\!
   \inttU e^{-t\gl}\,\go_1(r)k_a(\gl,r,x,r,x)\dbar\gl dr.
\end{equation*}
Since $r\ge t^{1/m}$ and $|\gl|\ge t^{-1}$, we may assume $r^m|\gl|\ge 1$.

\begin{lemma}\label{tr-as5}
For every $x\in X$, $0<t<1$ and $N\in\N$ with $N>j$, we get
\[ s''[a](t,x)=\sum_{k=0}^{N+n} c_{k}(x)\,t^{(k-n-1)/m}
   +\sum_{k=0}^{N-1} c_{k}'(x)\,t^{k/m}\log t+O(t^{N/m}), \]
where $n=\dim X$. The coefficients $c_k$ and $c_k'$ depending also on
$\ell$ and $j$.
\end{lemma}
\begin{proof}
The proof of this lemma is rather long. For this reason, we wish to
outline our strategy in order to make the structure of the proof clearer.
First of all, let $N>j$ be given, and recall that $\ell$ and $j$ are fixed. 
Let us simplify the notation by writing
\[ S(t,x)=s''[a](t,x)\;\text{ and }\;
   K(\gl,r,x)=r^{-m\ell-j}k_a(\gl,r,x,r,x). \]
In local coordinates the function $K(\gl,r,x)$ is given by
\[ K(\gl,r,x)=r^{-1}\int p(r,x,\eta,\gl)\dbar\eta \]
where $\eta=(\xi,z)\in \R^n\times\gG_\gb \cong\R^{n+1}$, 
$\dbar\eta=\frac{1}{(2\pi)^{n+1}}d\eta$, 
$p(r,x,\eta,\gl)=\tilde p(x,\eta,r^m\gl)$, and 
$\tilde p\in S_{\mathcal{O},\cl}^{-m\ell,m}(\gO\times\R^n\times\C;\gL)$
is a local symbol of $\tilde h(z,\gl)$. Thus there are symbols $\tilde p_k$,
anisotropic homogeneous in $(\eta,\gl)$ of degree $(-m\ell-k,m)$ for
$|\eta|\ge 1$, such that for every $J\in\N$,
\[ \tilde p=\sum_{k=0}^{J-1}\tilde p_k + \tilde g_J
   \;\text{ with } \tilde g_J
   \in S_{\mathcal{O}}^{-m\ell-J,m}(\gO\times\R^n\times\C;\gL).\]
This expansion induces a decomposition of $K(\gl,r,x)$, say
$K=K_{p_0}+K_{p_1}+\cdots+K_{g_J}$, and a decomposition of $S(t,x)$,
say $S=S_{p_0}+\cdots +\,S_{p_{J-1}}+\,S_{g_J}$. Our purpose is to expand
in $t$ all these components. 

\bigskip
{\it\underline{Step 1:}}\hspace*{1ex}(\textit{Expansion of} $S_{g_J}$). 
First of all, let $g_J(r,x,\eta,\gl)=\tilde g_J(x,\eta,r^m\gl)$. 
Following the proof of Lemma~\ref{tr-as3} we can write, 
exactly as in \eqref{tr-as3.1}, 
\begin{align*} 
S_{g_J}(t,x)
  &= t^{j/m}\left\{\int_1^{\infty}\!\!
  \intU e^{-\gl} k_J(\gl,r,x,r,x)\dbar\gl dr\right\}\\ \notag
  &\quad + t^{N/m}\int_1^{\infty}\!\!
  \intU e^{-\gl}\go_{1,N-j}(t^{1/m}r)\,r^{N-j} k_J(\gl,r,x,r,x)\dbar\gl dr,
\end{align*}
where $k_J=r^{m\ell+j}K_{g_J}(\gl,r,x)$. For $J>N$ the kernel $k_J$ decays 
fast enough as $r\to\infty$ so that the integrals exist, and the last one 
is uniformly bounded. Thus
\begin{equation} \label{tas4}
S_{g_J}(t,x) = c(x) t^{j/m} + O(t^{N/m}) 
\end{equation}
with $c(x)$ depending on $\ell, j, J$.

\bigskip
{\it\underline{Step 2:}}\hspace*{1ex}(\textit{Splitting of} $S_{p_k}$).
For $k<J$ let us split (recall that $r^m|\gl|\geq 1$)
\begin{equation}\label{tas5}
  \int p_k\,\dbar\eta=\int_{|\eta|\geq r|\gl|^{1/m}}p_k\,\dbar\eta+
  \int_{|\eta|\leq 1}p_k\,\dbar\eta+
  \int_{r|\gl|^{1/m}\geq|\eta|\geq 1}p_k\,\dbar\eta.
\end{equation}
This equality split likewise $S_{p_k}=S_{k,1}+S_{k,2}+S_{k,3}$ with
\begin{equation}\label{tas5.1}
S_{k,\gd} = t^{-\ell+1}\!\int_{t^{1/m}}^\infty\!    
 \inttU e^{-t\gl} \go_1(r) r^{m\ell+j-1}
 \left(\int_{V_\gd}\tilde p_k(r,x,\eta,r^m\gl)\dbar\eta\right)\dbar\gl dr
\end{equation}
for $\gd\in\{1,2,3\}$, where $V_\gd$ denotes the corresponding domain of 
integration. 

\bigskip
{\it\underline{Step 3:}}\hspace*{1ex}(\textit{Expansion of} $S_{k,1}$). 
The change of variables $\eta\to r|\gl|^{1/m}\eta\,$ yields 
\begin{align*}
 \int_{|\eta|\geq r|\gl|^{1/m}}\tilde p_k(x,\eta,r^m\gl)\dbar\eta
   &=\int_{|\eta|\geq 1}(r|\gl|^{1/m})^{n+1}
     \tilde p_k(x,r|\gl|^{1/m}\eta,r^m|\gl|\tfrac{\gl}{|\gl|})\dbar\eta\\
   &=(r|\gl|^{1/m})^{-m\ell-k+n+1}\int_{|\eta|\geq 1}
     \tilde p_k(x,\eta,\tfrac{\gl}{|\gl|})\dbar\eta
\end{align*}
due to the homogeneity of $\tilde p_k$. If we denote 
$\tilde c_k(x,\gl/|\gl|) =\int_{|\eta|\geq 1}
\tilde p_k(x,\eta,\gl/|\gl|) \dbar\eta$, and set $m_k:=-m\ell-k+n+1$, 
then for $\gb=1$ the integral \eqref{tas5.1} becomes
\begin{equation}\label{tas6.0}
\begin{split}
S_{k,1} &= t^{-\ell+1}\!\int_{t^{1/m}}^\infty\!    
 \inttU e^{-t\gl} \go_1(r) r^{m\ell+j-1} \left(
 (r|\gl|^{1/m})^{m_k} \tilde c_k(x,\gl/|\gl|)\right) \dbar\gl dr \\
 &= t^{(k-n-1)/m}\!\int_{t^{1/m}}^\infty\!   
 \intU e^{-\gl} \go_1(r) r^{j-k+n} |\gl|^{m_k/m} 
 \tilde c_k(x,\gl/|\gl|) \dbar\gl dr \\
 &= t^{(k-n-1)/m}\!\int_{t^{1/m}}^\infty\! \go_1(r) r^{j-k+n} dr 
 \left(\intU e^{-\gl} |\gl|^{m_k/m}\tilde c_k(x,\gl/|\gl|)\dbar\gl \right). 
\end{split}
\end{equation}
Finally, Lemma~\ref{tr-as6} leads to the expansion
\begin{equation}\label{tas6}
S_{k,1}(t,x)=\bar{c}(x)t^{(k-n-1)/m}+c(x)t^{j/m}+c'(x)t^{j/m}\log t
\end{equation}
with coefficients depending on $\ell,j,k$. In particular, $c'=0$
for $k\not=j+n+1$.

\bigskip
{\it\underline{Step 4:}}\hspace*{1ex}(\textit{Expansion of} $S_{k,2}$). 
Since $\tilde p_k\in S_{\mathcal{O}}^{-m\ell-k,m}(\gO\times\R^n\times\C;\gL)$
depends holomorphically on $\gl$, we can apply Theorem~\ref{p-hs7}
with $d=m$, $\nu=-k$ and $\mu=m\ell$ in order to expand $\tilde p_k$ as 
\[ \tilde p_k(x,\eta,\gl)=
   \sum_{\ga=0}^{M-1}\gl^{-\ell-\ga/m} q_{k\ga}(x,\eta)
   +\gl^{-\ell-M/m}\tilde p_{k,M}(x,\eta,\gl) \]
for any $M\in\N$, with $q_{k\ga}\in S_{\mathcal{O}}^{-k+\ga}
(\gO\times\R^n\times\C)$ and $\tilde p_{k,M}(\cdot,\cdot,\gl)\in 
S_{\mathcal{O}}^{-k+M}(\gO\times\R^n\times\C)$ uniformly for $\gl\in\gL$ with
$|\gl|\geq 1$. For $r^m|\gl|\ge 1$, the degenerate symbol 
$p_k(r,x,\eta,\gl)=\tilde p_k(x,\eta,r^m\gl)$ can then be written as 
\begin{equation*}
 p_k= \sum_{\ga=0}^{M-1}(r^m\gl)^{-\ell-\ga/m} q_{k\ga}(x,\eta)
     +(r^m\gl)^{-\ell-M/m}\tilde p_{k,M}(x,\eta,r^m\gl), 
\end{equation*}
and so
\begin{equation}\label{tas7}
\begin{split}
 \int_{V_2}p_k\dbar\eta =& \sum_{\ga=0}^{M-1} (r^m\gl)^{-\ell-\ga/m}
 \int_{V_2}q_{k\ga}(x,\eta)\dbar\eta \\
 &+ (r^m\gl)^{-\ell-M/m} \int_{V_2}\tilde p_{k,M}(x,\eta,r^m\gl)\dbar\eta.
\end{split}
\end{equation}
Set $\tilde c_{k,M}(x,r^m\gl)=\int_{V_2}\tilde p_{k,M}(x,\eta,r^m\gl)
\dbar\eta\,$ which is finite as $V_2=\{|\eta|\le 1\}$ is compact.
Now, since $\intU e^{-\gl}\gl^{-\ell-\ga/m}\dbar\gl=0$, we get 
\begin{equation}\label{tas3.0}
\begin{split}
S_{k,2}
 &= t^{M/m}\!\int_{t^{1/m}}^\infty\! \intU e^{-\gl} \go_1(r) r^{j-M-1} 
 \gl^{-\ell-M/m} \tilde c_{k,M}(x,r^mt^{-1}\gl) \dbar\gl dr\\
 &= t^{j/m}\!\int_{1}^\infty\! \intU e^{-\gl} \go_1(t^{1/m}r) r^{j-M-1} 
 \gl^{-\ell-M/m} \tilde c_{k,M}(x,r^m\gl) \dbar\gl dr
\end{split}
\end{equation}
due to the change of variables $\gl\to t^{-1}\gl$ and $r\to t^{1/m}r$.
Similarly as in the proof of Lemma~\ref{tr-as3} we expand
$\go_1(r)=1+r^{M-j-1}\go_{1,\tilde M}(r)$. Thus
\begin{align*}
S_{k,2}
 &= t^{j/m}\!\int_{1}^\infty\! \intU e^{-\gl} r^{j-M-1}
 \gl^{-\ell-M/m} \tilde c_{k,M}(x,r^m\gl) \dbar\gl dr\\
 &\quad + t^{(M-1)/m}\!\int_{1}^\infty\! r^{-2}\!
 \intU e^{-\gl} \go_{1,\tilde M}(t^{1/m}r)
 \gl^{-\ell-M/m} \tilde c_{k,M}(x,r^m\gl) \dbar\gl dr
\end{align*}
Since $\go_{1,\tilde M}(t^{1/m}r)$ and $\tilde c_{k,M}(x,r^m\gl)$ are
uniformly bounded for $0<t<1$, $r\geq 1$ and $r^m|\gl|\geq 1$, then the
integral at $t^{(M-1)/m}$ exists (presence of $r^{-2}$) and is also uniformly
bounded. Choosing $M\geq N+1$ we finally get  
\begin{equation}\label{tas8}
 S_{k,2}(t,x)=c(x)t^{j/m} + O(t^{N/m})
\end{equation}
with $c(x)$ depending on $\ell, j, k$.

\bigskip
{\it\underline{Step 5:}}\hspace*{1ex}(\textit{Expansion of} $S_{k,3}$). 
According to \eqref{tas7} let
\begin{equation*}
\begin{split}
 I_{\ga} &=(r^m\gl)^{-\ell-\ga/m}\int_{V_3}q_{k\ga}(x,\eta)\dbar\eta,\\
\tilde I_{M} 
   &=(r^m\gl)^{-\ell-M/m}\int_{V_3}\tilde p_{k,M}(x,\eta,r^m\gl)\dbar\eta,
\end{split}
\end{equation*}
so that
\[ S_{k,3}=t^{-\ell+1}\!\int_{t^{1/m}}^\infty\!\inttU \etgl\go_1(r) 
   r^{m\ell+j-1} \left(\sum_{\ga=0}^{M-1}I_{\ga} +
   \tilde I_{M}\right)\dbar\gl dr. \]
Recall that in \eqref{tas7} each $q_{k\ga}(x,\eta)=\frac{1}{\ga!}
\lim_{w\to 0}\partial_w^\ga\bigl(w^{-m\ell}
\tilde p_k(x,\eta,\frac{1}{w^m})\bigr)$ is homogeneous in $\eta$ of degree
$-k+\ga$ and consequently $\tilde p_{k,M}$ is anisotropic homogeneous of 
degree $(-k+M,m)$ for $|\eta|\geq 1$. Then
\[ \int_{V_3}\tilde p_{k,M}(x,\eta,r^m\gl)\dbar\eta
  =\int_{V_3} |\eta|^{-k+M}
  \tilde p_{k,M}(x,\tfrac{\eta}{|\eta|},r^m\tfrac{\gl}{|\eta|^m})\dbar\eta. \]
Taking now $M$ such that $-k+M>-(n+1)$ the latter integral can be split 
in two integrals $\int_{r|\gl|^{1/m}\geq|\eta|}-\int_{|\eta|\leq 1}$.
We denote the second integral again by $\tilde c_{k,M}(x,r^m\gl)$, and make 
the change of variables $\eta\to r|\gl|^{1/m}\eta$ in the first one. Thus
\begin{align*}
 \int_{r|\gl|^{1/m}\ge |\eta|}\tilde p_{k,M}\dbar\eta
 &= (r|\gl|^{1/m})^{-k+M+n+1}\int_{|\eta|\le 1}|\eta|^{-k+M}
 \tilde p_{k,M}(x,\tfrac{\eta}{|\eta|},\tfrac{\gl}{|\gl||\eta|^m})\dbar\eta\\ 
 &= (r|\gl|^{1/m})^{-k+M+n+1}\tilde c_{k,M}(x,\gl/|\gl|).
\end{align*}
Now, with the usual change of variables of $\gl$ and $r$,
\begin{multline}\label{tas9}
t^{-\ell+1}\!\int_{t^{1/m}}^\infty\!\inttU \etgl\go_1(r)
   r^{m\ell+j-1} \tilde I_{M}(r,x,\gl)\dbar\gl dr\\ 
\shoveright{
 =t^{(k-n-1)/m}\!\int_{t^{1/m}}^\infty\!\go_1(r)r^{j-k+n} dr
 \left(\intU e^{-\gl} |\gl|^{m_k/m}
 \tilde{\tilde c}_{k,M}(x,\tfrac{\gl}{|\gl|})\dbar\gl\right)} \\ 
 -t^{j/m}\!\int_{1}^\infty\!\intU e^{-\gl}\go_1(t^{1/m}r)r^{j-M-1}
  \gl^{-\ell-M/m} \tilde c_{k,M}(x,r^m\gl)\dbar\gl dr
\end{multline}
with $m_k=-m\ell-k+n+1$ and $\tilde{\tilde c}_{k,M}=
(\gl/|\gl|)^{-\ell-M/m} \tilde c_{k,M}(x,\gl/|\gl|)$.
These integrals are of the same type as \eqref{tas6.0} and \eqref{tas3.0},
so \eqref{tas9} is of the form
\begin{equation}\label{tas9.1} 
  \bar{c}(x)t^{(k-n-1)/m}+c(x)t^{j/m}+c'(x)t^{j/m}\log{t}+O(t^{N/m})
\end{equation}
with coefficients depending on $\ell,j,k$. 
In particular, $c'=0$ for $k\not=j+n+1$.

In order to complete the expansion of $S_{k,3}$ we next treat those terms
involving $I_{\ga}$. By means of polar coordinates and the homogeneity 
in $\eta$ of $q_{k\ga}$, 
\begin{multline*}
\int_{V_3}q_{k\ga}(x,\eta)\dbar\eta
  =d_{k\ga}(x)\int_1^{r|\gl|^{1/m}}\!\!\gr^{-k+\ga+n}d\gr\\ 
  =d_{k\ga}(x)
  \begin{cases}
   \tfrac{1}{\ga+n+1-k}\left((r|\gl|^{1/m})^{-k+\ga+n+1}-1\right) 
                       & \text{for } \ga\not=k-n-1\\
   \log(r|\gl|^{1/m})  & \text{for } \ga=k-n-1
  \end{cases}
\end{multline*}
which we denote by $\tilde c_{k\ga}(x,r|\gl|^{1/m})$.
Then, for $\ga<M$, the corresponding component of $S_{k,3}$ becomes
\begin{align} \notag
 t^{-\ell+1}\!\int_{t^{1/m}}^\infty\!\inttU 
 & \etgl\go_1(r) r^{m\ell+j-1} I_{\ga}(r,x,\gl)\dbar\gl dr \\  \label{tas10}
 &= t^{\ga/m}\!\int_{t^{1/m}}^\infty\!\intU
   e^{-\gl}\go_1(r) r^{j-\ga-1} \gl^{-\ell-\ga/m} 
   \tilde c_{k\ga}(x,r|t^{-1}\gl|^{1/m})\dbar\gl dr.
\end{align}
If $\ga\not=k-n-1$, then \eqref{tas10} equals
\[ t^{(k-n-1)/m}\int_{t^{1/m}}^\infty\!\go_1(r)r^{j-k+n}dr 
   \left(\tfrac{d_{k\ga}(x)}{\ga+n+1-k}
   \intU e^{-\gl}\gl^{-\ell-\ga/m}|\gl|^{(\ga+n+1-k)/m}\dbar\gl\right) \]
since $\intU e^{-\gl}\gl^{-\ell-\ga/m}\dbar\gl=0$ for every $\ga$.
Making use of Lemma~\ref{tr-as6} this expression can be written
as a linear combination of $t^{(k-n-1)/m}$, $t^{j/m}$ and $t^{j/m}\log t$.
Moreover, the term $t^{j/m}\log t$ disappears when $k\not=j+n+1$.
If $\ga=k-n-1$, then \eqref{tas10} becomes
\begin{equation*}
   t^{(k-n-1)/m}\int_{t^{1/m}}^\infty\!\go_1(r)r^{j-k+n}dr
   \left(d_{k\ga}(x) \intU e^{-\gl}\gl^{-\ell-\ga/m}
   \log|\gl|^{1/m}\dbar\gl\right) 
\end{equation*}
since the integral $\int_{t^{1/m}}^\infty\!\go_1(r)r^{j-k+n}\log(t^{-1/m}r)dr$
is finite and $\intU e^{-\gl}\gl^{-\ell-\ga/m}\dbar\gl=0$. Applying 
Lemma~\ref{tr-as6} once more we finally obtain
\begin{equation}\label{tas13}
 S_{k,3}(t,x)= 
   \bar{c}(x)t^{(k-n-1)/m}+c(x)t^{j/m}+c'(x)t^{j/m}\log t+O(t^{N/m}),
\end{equation}
with coefficients depending on $\ell,j,k$. 
In particular, $c'=0$ for $k\not=j+n+1$.

\bigskip
{\it\underline{Summary:}}\hspace*{1ex} 
By virtue of the relation
\[ s''[a](t,x)=\sumk{J-1}\bigl(S_{k,1}(t,x)+S_{k,2}(t,x)+S_{k,3}(t,x)\bigr)
          + S_{g_J}(t,x), \]
the assertion follows summing up \eqref{tas6}, \eqref{tas8}, \eqref{tas13}
and \eqref{tas4}. 
\end{proof}

\begin{lemma}\label{tr-as6}
Let $0<\tau<1$ and let $\go$ be a cut-off function with $\go(r)=1$ for $r<1$.
For every $j$, $\nu\in\R$ there is a constant $C_{j,\nu}>0$ such that
\[ \tau^\nu\int_\tau^\infty \go(r)r^{j-\nu-1}dr= C_{j,\nu}\,\tau^\nu +
    \begin{cases}
       -\tfrac{1}{j-\nu}\,\tau^j & \text{for } \nu\not=j\\
       -\tau^j\log \tau & \text{for } \nu=j
    \end{cases} \]
\end{lemma}
\begin{proof}
Since $\go(r)=1$ on $[0,1)$,
\begin{align*}
  \tau^\nu\int_\tau^\infty \go(r)r^{j-\nu-1}dr
   &=\tau^\nu\left(\int_1^\infty \go(r)r^{j-\nu-1}dr
       +\int_\tau^{1} r^{j-\nu-1}dr\right)\\
   &=\tau^\nu\int_1^\infty \go(r)r^{j-\nu-1}dr + 
    \begin{cases}
      \tfrac{\tau^\nu}{j-\nu}-
      \tfrac{\tau^j}{j-\nu} & \text{for } \nu\not=j\\
      -\tau^j\log \tau      & \text{for } \nu=j 
    \end{cases}
\end{align*} 
\end{proof}

In order to analyze the contribution of $R_N(t)$ from \eqref{R_N}, observe
that the Mellin operator involved in its definition can be written as 
$\,a_N(\gl)=r^{m\ell+N}\op_M(h_N)(\gl)$ with
\begin{equation*} 
 h_N(r,z,\gl) = \tfrac{1}{(N-1)!} \int_0^1 (1-\gt)^{N-1} 
 \tilde h(\gt r,z,r^m\gl)d\gt,
\end{equation*}
where $\tilde h\in C^\infty(\R_+,M_{\mathcal{O}}^{-m\ell,m}(X;\gL))$.
As seen in the previous calculations, there is no loss of generality if 
we drop the variable $x\in X$. Let $k_{a_N}(\gl,r,r')$ be the continuous 
Schwartz kernel of $a_N(\gl)$ and let
\begin{equation*}
 s[a_N](t) = t^{-\ell+1}\!\int_{0}^\infty\!
 \inttU e^{-t\gl}\,\go_1(r)k_{a_N}(\gl,r,r)\dbar\gl dr.
\end{equation*}
Split $s[a_N](t)=s'[a_N](t)+s''[a_N](t)$ as done before. As a matter of fact,
the proof of Lemma~\ref{tr-as5} also works for $s''[a_N](t)$ and we get 
\begin{equation*}
 s''[a_N](t)=\sum_{k=0}^{N+n} c_{k}\,t^{(k-n-1)/m} + O(t^{N/m}). 
\end{equation*}
The main difference between $a_N(\gl)$ and the Mellin operators considered
 at the beginning of this section, is the additional dependence on $r$ of 
the nondegenerate symbol $\tilde h(r,z,\gl)$. The family $a_N(\gl)$ is 
unfortunately not twisted homogeneous and Lemma~\ref{tr-as4} 
cannot be applied. However, we have the following 

\begin{lemma} \label{tr-as7}
Let $0<t<1$ and $N\in\N$. Then, $s'[a_N](t)= O(t^{N/m})$.
\end{lemma}

\begin{proof}
We can write $k_{a_N}(\gl,r,r)= r^{m\ell+N} K_N(\gl,r)$ with
\begin{equation*}
 K_N(\gl,r) = \tfrac1{(N-1)!}\, r^{-1}\! \iint_0^1 (1-\gt)^{N-1}
 \tilde h(\gt r,\tau,r^m\gl)d\gt \dbar\tau.
\end{equation*}
Since $(1-\gt)^{N-1} = \sum_{k=0}^{N-1} \binom{N-1}{k}(-1)^k \gt^k$, 
the function $K_N(\gl,r)$ can be split, inducing a corresponding decomposition 
$s'[a_N](t) =\sum_{k=0}^{N-1} c_{N,k} S_k(t)$ with
\begin{align*}
 S_k(t) &= t^{-\ell+1}\!\int_0^{t^{1/m}}\!\!
 \inttU e^{-t\gl}\,r^{m\ell+N-1}\! \left(\iint_0^1 \gt^k
 \tilde h(\gt r,\tau,r^m\gl)d\gt \dbar\tau\right) \dbar\gl dr \\
 &= t^{(N-1-k)/m}\!\int_0^{1}\!
 \intU e^{-\gl}\,r^{m\ell+N-1}\! \left(\iint_0^{t^{1/m}} \!\gt^k
 \tilde h(\gt r,\tau,r^m\gl)d\gt \dbar\tau\right) \dbar\gl dr 
\end{align*}
after the change of variables $\gl\to t^{-1}\gl$, $r\to t^{1/m}r$, 
and $\gt\to t^{-1/m}\gt$. Now, 
\begin{equation*}
 \int_0^{t^{1/m}}\! \gt^k \tilde h(\gt r,\tau,r^m\gl)d\gt 
   = \tfrac1{k+1} t^{(k+1)/m} \tilde h(\xi_t r,\tau,r^m\gl)
\end{equation*}
for some $\xi_t\in [0,t^{1/m}]$, and so  
\begin{equation*}
 S_k(t) = t^{N/m} \tfrac1{k+1}\int_0^{1}\!
 \intU e^{-\gl}\,r^{m\ell+N-1}\! \left(\int
 \tilde h(\xi_t r,\tau,r^m\gl) \dbar\tau\right) \dbar\gl dr 
\end{equation*}
which is $O(t^{N/m})$ since the integral is uniformly bounded.
\end{proof}

\begin{theorem}\label{htrace-exp}
Let $A$ be a cone differential operator of order $m>0$ such that $A-\gl$ is 
parameter-elliptic with respect to some $\gg\in\R$ on a sector $\gL$ as given
in Section~\ref{heat-op}. Then, the heat trace admits the asymptotic expansion
\[ \quad\tr e^{-tA}\sim \sum_{k=0}^\infty C_k\,t^{(k-n-1)/m}+
   \sum_{k=0}^\infty C_k'\,t^{k/m}\log{t}, \;\text{ as }\; t\to 0^+, \]
where $C_k$ and $C_k'$ are constants depending on the symbolic structure of $A$.
\end{theorem}

\begin{proof}
Actually, we only need to put together all the single expansions obtained in
this section. The asymptotic summation here means:

\smallskip
For a given $K\in\N$ there exists $N(K)\in\N$ and $C_K>0$ such that
\[ \Big|\tr\etA-\sumk{N-1}\tau_k(t)\Big|\leq C_K\,t^{K} \]
for every $N\geq N(K)$ and $0<t<1$, where $\tau_k(t)$ is some expression 
of 'degree' $k$ from the right-hand side above, that is, $\tau_k(t)$ is
a linear combination of terms like $t^{(k-\ga)/m}\log^\gb{t}$ with
$\ga=0\text{ or } n+1$, and $\gb=0\text{ or } 1$. 

First of all, we choose $\ell\in\N$ large enough such that $(A-\gl)^{-\ell}$
is of trace class. Hence the family $E_\ell^{(N)}(t)$ from \eqref{h-e1} is of
trace class too, and Corollary~\ref{h-as3} together with Theorem~\ref{trace}
 imply: for a given $K\in\N$ there are $N(K)$ and $C_K>0$ such that
\[ \left|\tr\etA - \tr E_\ell^{(N)}(t)\right|\le C_K\,t^{K} \]
for $0<t<1$. Recall that $E_\ell^{(N)}(t)=\sumk{N-1} E_{\ell,k}(t)+ R_N(t)$
with $E_{\ell,k}(t)$ and $R_N(t)$ as in \eqref{R_N}. 
By definition, every $E_{\ell,k}(t)$ splits into a component supported in the 
interior of $\mathbb{B}$ and another one supported near $\partial\mathbb{B}$.
It is known (cf. \cite[Sec. 1.8.1]{Gilkey}) that the trace of the interior 
part is of the desired form due to the homogeneity of the local symbols.
On the other hand, all the components supported near the boundary are either  
operator-valued Green symbols or parameter-dependent Mellin operators with 
degenerate symbols. By using Lemma~\ref{tr-as3} for the integrals of Green 
type and Lemma~\ref{tr-as4}, Lemma~\ref{tr-as5} and Lemma~\ref{tr-as7} for 
the rest, we obtain the desired expansion integrating 
over the manifold $X$ the local expansions obtained there.
\end{proof}

\renewcommand\thesubsection{\Alph{subsection}}
\renewcommand\thetheorem{\Alph{subsection}.\arabic{theorem}}
\section{Appendix}
\setcounter{subsection}{0}
\subsection{Parameter-dependent pseudodifferential operators} \label{pdo-p}
In this section we give the basic definitions of the local parameter-dependent
symbols that we use in this paper. For further information the reader is
referred, for instance, to the book of {\sc Shubin} \cite[Section~9]{Shubin}.

\smallskip
Let $\gO$ be an open set in $\subset\R^n$, and let $\gL$ be a closed angle 
in $\C$ with vertex at the origin.
Let $\mu\in\R$, $d\in\N$. A function $p(x,\xi,\gl)\in
C^{\infty}(\gO\times\R^n\times\gL)$ is said to be in the class of symbols 
$S^{\mu,d}(\gO\times\R^n;\gL)$ if for any multi-indices 
$\ga,\gb\in\N_0^{n}$, $\gg\in\N_0^2$, and any compact set $K\subset\gO$ 
there is a positive constant $C_{\ga,\gb,\gg,K}$  such that 
\[ |D^\ga_\xi D^\gb_x D^\gg_\gl p(x,\xi,\gl)|\leq 
   C_{\ga,\gb,\gg,K}\big(1+|\xi|+|\gl|^{1/d}\big)^{\mu-|\ga|-d|\gg|} \]
for $x\in K$, $\xi\in\R^n$, $\gl\in\gL$. The smoothing elements
\[ S^{-\infty}(\gO\times\R^n;\gL):=
   \bigcap_{\mu\in\R}S^{\mu,d}(\gO\times\R^n;\gL)\]
are independent of $d$. To every symbol $p\in S^{\mu,d}(\gO\times\R^n;\gL)$
we associate the operator family $P(\gl)=\op_x(p)(\gl)$ given by
\[ [P(\gl)u](x):=\iint e^{i(x-y)\xi}p(x,\xi,\gl)u(y)\,dy\,\dbar\xi
   \;\text{ for } u\in C^{\infty}_{0}(\gO). \]
It is a pseudodifferential operator of order $\mu$ depending smoothly on the
parameter $\gl\in\gL$ with  anisotropy $d$; we write $P(\gl)\in
L^{\mu,d}(\gO;\gL)$. When $d=1$ we omit it from the notation.

\smallskip
A function $f(x,\xi,\gl)$ on $\gO\times\R^n\times\gL$ is called
{\em (anisotropic) homogeneous} in $(\xi,\gl)$ of degree $(\nu,d)$ if
\[ f(x,\tau\xi,\tau^d\gl)=\tau^{\nu}\,f(x,\xi,\gl)\;
   \text{ for every } \tau\ge 1. \]

The class $S_{\cl}^{\mu,d}(\gO\times\R^n;\gL)$ of {\em classical} 
pseudodifferential symbols is defined as the space of symbols 
$p\in S^{\mu,d}(\gO\times\R^n;\gL)$ that admit an asymptotic expansion 
of the form
\[ p(x,\xi,\gl)\sim \sum_{j=0}^\infty p_{\mu-j}(x,\xi,\gl),\]
where each $p_{\mu-j}$ is homogeneous of degree $(\mu-j,d)$ 
for $|\xi|+|\gl|^{1/d}\ge 1$. 

\smallskip
For $\bar d=(d_1,d_2)\in\N^2$ a smooth function $p$ is said to be in
$S^{\mu,\bar d}(\gO\times\R^n;\R^\ell\times\gL)$, if for any
$\ga,\gb\in\N_0^{n}$, $\gg\in\N_0^{\ell+1}$, and any compact $K\subset\gO$
there is a positive constant $C=C({\ga,\gb,\gg,K})$ such that
\[ |D^\ga_\xi D^\gb_x D^{\gg_1}_\gr D^{\gg_2}_\gl p(x,\xi,\gr,\gl)|
   \leq C\big(1+|\xi|+|\gr|^{1/d_1}+|\gl|^{1/d_2}\big)^{\mu-|\ga|
   -d_1|\gg_1|-d_2\gg_2} \]
for $x\in K$, $\xi\in\R^n$ and $(\gr,\gl)\in\R^\ell\times\gL$. 

\subsection{Green cone operators}\label{s-cone}
\setcounter{theorem}{0}
In this section we only discuss a special class of Green operators that
describes sufficiently well the smoothing elements appearing in our context.
For recent and more general results concerning Green cone operators we
refer to \cite{Krainer}.

Let us begin by pointing out some embedding properties of the weighted Sobolev
spaces on the manifold $\mathbb{B}$ (cf. Section~\ref{sobolev}).   

\begin{lemma}\label{emb2} 
For $s\geq s',\gg\geq\gg'$ the embedding
$\mathcal{H}^{s,\gg}(\mathbb{B})\emb\mathcal{H}^{s',\gg'}(\mathbb{B})$
is continuous. If $\gg>\gg'$, then it is compact if $s>s'$, Hilbert-Schmidt 
if $s>s'+\dim\mathbb{B}/{2}$, and trace class when $s>s'+\dim\mathbb{B}$.
\end{lemma}

The intersection over $s\in\R$ of all these spaces will be denoted by 
$\mathcal{H}^{\infty,\gg}(\mathbb{B})$, and will be topologized as the
projective limit $\mathcal{H}^{\infty,\gg}(\mathbb{B})
 =\projlim_{s\in\N}\hsg(\mathbb{B})$ for every $\gg\in\R$. 
Finally, we introduce the space
\begin{equation}\label{hinf-m}
  \mathcal{H}^{\infty,\gg^-}(\mathbb{B})
  :=\projlim_{k\in\N}\mathcal{H}^{k,\gg-\frac{1}{k}}(\mathbb{B}). 
\end{equation}
Observe that $\mathcal{H}^{\infty,\gg}(\mathbb{B})
\emb\mathcal{H}^{\infty,\gg^-}(\mathbb{B})$. Moreover, 
$\mathcal{H}^{\infty,\gg^-}(\mathbb{B})$ is nuclear since the embedding
$\mathcal{H}^{k',\gg-\frac{1}{k'}}(\mathbb{B})
\emb\mathcal{H}^{k,\gg-\frac{1}{k}}(\mathbb{B})$ is Hilbert-Schmidt for
$k'>k+\frac{(n+1)}{2}$.

\begin{definition}\label{green1} 
Let $\gg,\gd\in\R$. An operator $G\in\bigcap_{s,s'\in\R}\mathcal{L}
\big(\mathcal{H}^{s,\gg}(\mathbb{B}),\mathcal{H}^{s',\gd}(\mathbb{B})\big)$
is called a {\it Green operator} if there is an $\eps=\eps(G)>0$ such that
\[ G:\mathcal{H}^{s,\gg}(\mathbb{B})\to
   \mathcal{H}^{s',\gd+\eps}(\mathbb{B})\;\text{ and }\; 
   G^*:\mathcal{H}^{s,-\gd}(\mathbb{B})\to
   \mathcal{H}^{s',-\gg+\eps}(\mathbb{B}) \]
are continuous maps for all $s,s'\in\R$, where $G^*$ is the formal adjoint of
$G$  with respect to $(\cdot,\cdot)_{\mathcal{H}^{0,0}}$. The space of Green
operators with asymptotic data $(\gg,\gd)$ is denoted by $\coneg{\gg,\gd}$,
or $\coneg{\gg,\gd}_\eps$ if we fix the value of $\eps>0$. The latter  
is a Fr{\'e}chet space with the system of (semi)norms $\{|\cdot|_\Sigma,\; 
|\cdot|_\Sigma^* \;\text{ for } \Sigma=(s,s')\in\mathbb{Z}^2\}$ given by
\[ |G|_\Sigma =\norm{G}_{\mathcal{L}(\mathcal{H}^{s,\gg},
   \mathcal{H}^{s',\gd+\eps})},\quad
   |G|_\Sigma^* =\norm{G^*}_{\mathcal{L}(\mathcal{H}^{s,-\gd},
   \mathcal{H}^{s',-\gg+\eps})}. \]
\end{definition}

\begin{example}\label{green1.1}
If $\go_1$ and $\go_2$ are cut-off functions on $\mathbb{B}$, then 
\[ (1-\go_1) L^{-\infty}(\iBB)(1-\go_2)\subset\coneg{\gg,\gd}_\eps \]
for every $\gg,\gd\in\R$ and $\eps>0$.
\end{example}

\begin{example}\label{green1.2}
Let $K\in\mathcal{H}^{\infty,-\gg+\eps}(\mathbb{B})\pitensor
\mathcal{H}^{\infty,\gd+\eps}(\mathbb{B})$ for $\gg,\gd\in\R$ and $\eps>0$, 
where $\pitensor$ denotes the projective tensor product. 
The operator $G$ defined by 
\[ Gu(y)=\int_\mathbb{B} K(y,y')u(y')dy' \;
   \text{ for } u\in C^{\infty}_{0}(\mathbb{B}), \]
is an element of $\coneg{\gg,\gd}_\eps$. Here, $dy'$ denotes the measure 
induced by the fixed metric on $\mathbb{B}$.
\end{example}

\begin{theorem}\label{green2} 
Let $\gg,\gd\in\R$, $\eps>0$ and $G\in\coneg{\gg,\gd}_{\eps}$. For every
$s, s'\in\R$ the operator $G$ belongs to 
$\mathcal{L}_2(\mathcal{H}^{s,\gg}(\mathbb{B}),
\mathcal{H}^{s',\gd}(\mathbb{B}))$, the space of Hilbert-Schmidt operators.
Moreover, its kernel satisfies 
\[ K_G\in\mathcal{H}^{\infty,(-\gg+\eps)^-}(\mathbb{B})
   \pitensor \mathcal{H}^{\infty,(\gd+\eps)^-}(\mathbb{B}). \]
\end{theorem}

In fact, $G$ is in $\mathcal{L}_2(\mathcal{H}^{s,\gg},
\mathcal{H}^{s',\gd})$ because it maps 
$\mathcal{H}^{s,\gg}\to\mathcal{H}^{\infty,\gd+\eps}$ 
continuously, and the embedding $\mathcal{H}^{\infty,\gd+\eps}\emb
\mathcal{H}^{s',\gd}$ is in $\mathcal{L}_2$ for every $s'\in\R$, 
see Lemma~\ref{emb2}. Observe that because of the relation  
\[ \mathcal{H}^{\infty,(-\gg+\eps)^-}\!(\mathbb{B})
   \pitensor\mathcal{H}^{\infty,(\gd+\eps)^-}\!(\mathbb{B})
   \emb C^{\infty}(\iBB)\pitensor C^{\infty}(\iBB)\cong
   C^{\infty}(\iBB\times\iBB), \]
we clearly have $\coneg{\gg,\gd}_\eps\subset L^{-\infty}(\iBB)$ 
for all $\gg,\gd$ and $\eps$.

For our purposes it is convenient to describe the kernels of
Green operators with the help of the Hilbert tensor product.
If $H_1$ and $H_2$ are Hilbert spaces then the tensor product
$H_1\otimes_{\scriptscriptstyle{H}}H_2$ is defined as the space of finite 
dimensional operators $H_1'\to H_2$ endowed with the topology of 
$\mathcal{L}_2(H_1',H_2)$. Hence 
\[ \qquad H_1\htensor H_2=\mathcal{L}_2(H_1',H_2)\quad(\text{completion}). \]
This concept of tensor product extends to the so-called hilbertizable
locally convex spaces, i.e., spaces whose completion can be written as a 
reduced projective limit of Hilbert spaces. Nuclear spaces are hilbertizable.
For further details we refer to \cite{DF93}, \cite{Jar}, \cite{Man95a}.
Let us just point out some basic properties. If $E$ and $F$ are hilbertizable
spaces then $E\htensor F$ is hilbertizable. Furthermore,  
$E\pitensor F\emb E\htensor F\emb E\hat{\otimes}_{\eps}F$, so
if $E$ or $F$ is a nuclear space, then
\[ E\pitensor F\cong E\hat{\otimes}_{\eps}F\cong E\htensor F. \]
Furthermore, the Hilbert tensor product commutes with projective limits, e.g.,
\[ \mathcal{H}^{\infty,\gg^-}\htensor\,\mathcal{H}^{\infty,\gd^-}
  =\projlim_{k,k'\in\N}\mathcal{H}^{k,\gg-\frac{1}{k}}
  \htensor\,\mathcal{H}^{k',\gd-\frac{1}{k'}}\;\text{ for } \gg,\gd\in\R.\]
According to Theorem~\ref{green2} the kernel of $G\in\coneg{\gg,\gd}_\eps$
belongs to
\[ \mathcal{H}^{\infty,(-\gg+\eps)^-}\pitensor
   \mathcal{H}^{\infty,(\gd+\eps)^-}\cong
   \mathcal{H}^{\infty,(-\gg+\eps)^-}\htensor
   \mathcal{H}^{\infty,(\gd+\eps)^-}
   \emb \mathcal{H}^{\infty,-\gg_\eps}\htensor\mathcal{H}^{\infty,\gd_\eps} 
\]
for every $\gd_\eps<\gd+\eps$ and $-\gg_\eps<-\gg+\eps$.

The embedding properties from Lemma~\ref{emb2} together with some standard 
results from functional analysis, cf. \cite[Appendix~3]{Shubin}, imply
the following
 
\begin{theorem}\label{trace}
Let $G$ be an operator of order $-\mu<-\dim\mathbb{B}$. Assume further that
$G:\hsg(\mathbb{B})\to \mathcal{H}^{s+\mu,\gg+\eps}(\mathbb{B})$ is bounded
for some $\gg\in\R$, $\eps>0$, and all $s\in\R$. Then $G$ is of trace class
on $\hsg(\mathbb{B})$ and, 
for any $k>s+\dim\mathbb{B}$ and $\gg\le \gg_\eps<\gg+\eps$, 
\[ |\tr G|\le C\norm{G}_{\mathcal{L}(\mathcal{H}^{-k,\gg},
   \mathcal{H}^{k,\gg_\eps})} \]
for some $C>0$. Notice that every $G\in\coneg{\gg,\gg}_\eps$ is of this form.
\end{theorem}

\subsection{A class of operator-valued symbols}\label{opv-symb}
\setcounter{theorem}{0}
One of the main ideas of {\sc Schulze}'s edge calculus is to consider the 
operators (near the edge) as operator-valued symbols acting on Banach
spaces together with groups of isomorphisms. In this section we introduce
analogously a class of anisotropic symbols similarly to \cite{BuSz97} but
here on conical subsets of the complex plane.
In this context, every Banach space $E$ will be considered together with
a {\it strongly continuous  group action} $\{\gk_{\tau}\}_{\tau\in\R_+}$, 
that is, a group of isomorphisms on $E$ satisfying
\begin{enumerate} 
 \item $\tau\mapsto\gk_{\tau}\,e:\R_+\to E$ is continuous for each $e\in E$,
 \item $\gk_{\tau}\gk_{\gs}=\gk_{\tau\gs}$ for all $\tau,\gs>0$.
\end{enumerate}

The following estimate can be obtained by means of the Banach-Steinhaus
theorem. For a proof see, e.g., Remark~2.2 in \cite{Hi90}.

\begin{lemma}\label{opvs0} 
There are constants $C,M\ge 0$ such that
\[ \norm{\gk_{\tau}}_{\mathcal{L}(E)}\leq C\max(\tau,\tau^{-1})^M. \]
\end{lemma}

Let us fix a positive smooth function $[\cdot]:\C \to \R_+$
with $[\gl]=|\gl|$ for $|\gl|\geq 1$. We also fix $d\in\N$ to describe the 
anisotropy. Further, let us set
\begin{equation}\label{opvs1.1}
 [\gl]_d:=[\gl]^{1/d}\;\text{ and }\; \gk(\gl):=\gk_{[\gl]_d} \;
 \text{ for some group action } \{\gk_\tau\}_{\tau\in\R_+}.  
\end{equation}

For the rest of this section let $\gL\subset\C$ be a conical set.
Let also $E_0$ and $E_1$ be Banach spaces with group actions 
$\gk_0=\{\gk_{0,\tau}\}_{\tau\in\R_+}$ and 
$\gk_1=\{\gk_{1,\tau}\}_{\tau\in\R_+}$, respectively.

\begin{definition}\label{opvs2} 
Let $\mu\in\R$. A function $a\in C^{\infty}(\gL,\mathcal{L}(E_0,E_1))$ is 
said to be in the space $S^{\mu,d}(\gL; E_0,E_1)$ of {\em operator-valued 
symbols}, if for every $\ga\in\N_0^2$ there is a constant $C_\ga>0$ such that
\begin{equation}\label{sdef}
  \norm{\gk_1^{-1}(\gl)\bigl\{\partial_{\gl}^\ga\, a(\gl)\bigr\} 
  \gk_0(\gl)}_{\mathcal{L}(E_0,E_1)}\leq C_\ga[\gl]_d^{\mu-d|\ga|}
\end{equation}
for all $\gl\in\gL$. Furthermore, we define
\[ S^{-\infty}(\gL; E_0,E_1):=\bigcap_{\mu\in\R} S^{\mu,d}(\gL; E_0,E_1)
   =\mathcal{S}(\gL,\mathcal{L}(E_0,E_1)) \]
which is independent of $d$ and the group actions $\gk_0$, $\gk_1$.
\end{definition}

Most of the usual symbol properties (e.g. embedding, composition, 
asymptotic summation, etc.) that are known for the scalar-valued symbols, 
can also be formulated in the operator-valued case (see e.g. \cite{BuSz97}).

\begin{proposition}\label{opvs3} 
Let $\mu,\tilde \mu\in \R$ and $E_2$, $E_3$ be further Banach spaces. Then
\begin{enumerate}
\item $S^{\mu,d}(\gL; E_1,E_2)\emb S^{\mu,d}(\gL; E_0,E_3)$ 
  if $E_0\emb E_1$ and $E_2\emb E_3$ with $\gk_1=\gk_0$ on $E_0$ and 
  $\gk_3=\gk_2$ on $E_2$,
\item $S^{\mu,d}(\gL;E_0,E_1)\emb 
  S^{\mu+M_0+M_1,d}(\gL;E_0,E_1)_{(\mathrm{id})}$,
  where $M_0$ and $M_1$ are associated to $\gk_0$ and $\gk_1$ as in 
  Lemma~\ref{opvs0}, and the subscript $(\mathrm{id})$ means that 
  $E_0$ and $E_1$ are considered with the trivial group action 
  $\kappa=\mathrm{id}$,
\item $\partial_\gl^\ga
S^{\mu,d}(\gL; E_0,E_1)\subset S^{\mu-d|\ga|,d}(\gL; E_0,E_1)$
  for every $\ga\in\N^2$,
\item $S^{\mu,d}(\gL; E_1,E_2) \cdot S^{\tilde \mu,d}(\gL; E_0,E_1)\subset
  S^{\mu+\tilde \mu,d}(\gL; E_0,E_2)$.
\end{enumerate}
\end{proposition}

An operator-valued function 
$a\in C^{\infty}(\gL,\mathcal{L}(E_0,E_1))$ is called
{\em twisted homogeneous of degree $(\mu,d)$} (with respect to the group
actions $\gk_0$ and $\gk_1$) if
\begin{equation}\label{t-hom} 
  a(\tau^d\gl) = \tau^\mu \gk_{1,\tau}\,a(\gl)\,\gk_{0,\tau}^{-1} 
\end{equation}
for all $\tau\ge 1$ and $\gl\in\gL$ sufficiently large. 

\begin{lemma}\label{tr-as2}
Let $E_0$ and $E_1$ be Banach spaces over $X^\wedge=\R_+\times X$ with the 
same group action. Let $a\in C^{\infty}(\gL,\mathcal{L}(E_0,E_1))$ be 
twisted homogeneous of degree $(\mu,d)$, and assume that every $a(\gl)$ is 
an integral operator with kernel $k_a(\gl,r,x,r',x')$. Then 
\begin{equation}\label{tr-a.e2}
 k_a(\tau^d\gl,r,x,r',x')=\tau^{\mu+1}k_a(\gl,\tau r,x,\tau r',x')
\end{equation}
for every $\tau\ge 1$ and $\gl\in\gL$ sufficiently large. 
\end{lemma}

The concept of twisted homogeneity allows to define classical symbols:

\begin{definition}\label{opvs5} 
The space $S_{\cl}^{\mu,d}(\gL;E_0,E_1)$ is the subspace of
$S^{\mu,d}(\gL;E_0,E_1)$ consisting of symbols $a(\gl)$ with an
expansion $a(\gl)\sim\sumj{\infty} a_{\mu-j}(\gl)$ such that every 
$a_{\mu-j}$ is twisted homogeneous of degree $(\mu-j,d)$.
In this case we call $\gs^\mu(a)(\gl):=a_\mu(\gl)$ the {\em principal part}
of $a(\gl)$.
\end{definition}

\begin{example}\label{opvs6} 
If $a\in C^{\infty}(\gL,\mathcal{L}(E_0,E_1))$ is twisted homogeneous
of degree $(\mu,d)$, then $a\in S_{\cl}^{\mu,d}(\gL;E_0,E_1)$.
\end{example}

\begin{proposition}\label{opvs7} 
Let the symbol $a\in S^{\mu,d}(\gL;E_0,E_1)$ be such that 
$a(\gl):E_0\to E_1$ is an invertible operator for every $\gl\in\gL$. Then 
$a^{-1}\in S^{-\mu,d}(\gL;E_1,E_0)$ if and only if 
$\norm{\gk_0^{-1}(\gl) a(\gl)^{-1}\gk_1(\gl)}_{\mathcal{L}(E_1,E_0)}
\le C[\gl]_d^{-\mu}$ for some $C>0$.
\end{proposition}

{\bf Green operator-valued symbols}

\medskip
Let $\mu,\gg,\gd\in\R$, $d\in\N$ and $\eps>0$. The function
\[ g\in\bigcap_{\Sigma\in\R^4}C^{\infty}
  \big(\gL,\mathcal{L}(\mathcal{K}^{s,\gg}(X^\wedge)^\gs,
  \mathcal{K}^{s',\gd}(X^\wedge)^{\gs'})\big), \]
where $\Sigma=(s,s',\gs,\gs')$, is said to be in the class 
$\edgeg{\mu,d}{\gg,\gd}_\eps\,$ if and only if
\begin{align*} 
 g &\in\bigcap_{\Sigma\in\R^4} S_{\cl}^{\mu,d}
  \big(\gL; \mathcal{K}^{s,\gg}(X^\wedge)^\gs, 
  \mathcal{K}^{s',\gd+\eps}(X^\wedge)^{\gs'}\big)\;\text{ and} \\
 g^* &\in\bigcap_{\Sigma\in\R^4} S_{\cl}^{\mu,d}
  \big(\gL; \mathcal{K}^{s,-\gd}(X^\wedge)^\gs, 
  \mathcal{K}^{s',-\gg+\eps}(X^\wedge)^{\gs'}\big),
\end{align*}
where $g^*(\gl)$ is the pointwise formal adjoint of $g(\gl)$ in 
$\mathcal{K}^{0,0}$.  Furthermore, we denote by $\edgeg{\mu,d}{\gg,\gd}$
 the union for $\eps>0$ of all these classes. Every element 
$g\in\edgeg{\mu,d}{\gg,\gd}$ is called a {\em Green symbol} of order
$(\mu,d)$ with weight data $(\gg,\gd)$. As a classical operator-valued 
symbol such a Green symbol $g$ has a twisted homogeneous principal component
which we denote $\sK^\mu(g)(\gl)$. Finally, we set
\[ \edgeg{-\infty}{\gg,\gd}:=\bigcap_{\mu\in\R} \edgeg{\mu,d}{\gg,\gd} \]
which is independent of $d$ and the group actions involved.

\section*{acknowledgements}
This paper is based on my dissertation \cite{Gil} which was submitted to  
the University of Potsdam in Germany. I would like to thank 
my advisors {\sc B.-W.~Schulze} and {\sc E.~Schrohe} for their excellent
guidance and constant support. The financial
support by the Max-Planck-Gesellschaft is greatly appreciated.
I thank {\sc T.~Krainer} and {\sc J.~Seiler} for many helpful 
discussions along the way. 
I also wish to thank {\sc M.~Lesch} for his help during my stay at 
the Humboldt University, and {\sc G.~Mendoza} for his constant
support at Temple University, and for our numerous and fruitful discussions.


\newcommand{\noopsort}[1]{}
\providecommand{\bysame}{\leavevmode\hbox to3em{\hrulefill}\thinspace}

\end{document}